\documentclass [reqno,12pt,a4paper]{amsart}
\usepackage[foot]{amsaddr}
\usepackage[a4paper,margin=1in,footskip=0.4in,heightrounded]{geometry}
\usepackage{lastpage}
\usepackage{t1enc}
\usepackage[utf8]{inputenc}
\usepackage{graphicx}
\usepackage{tabularx}
\usepackage{caption,subcaption}
\usepackage{multirow}
\usepackage{tikz,color}
\usepackage{mathtools}
\usepackage{latexsym,amssymb}
\usepackage{epstopdf,afterpage}
\usetikzlibrary{arrows}
\graphicspath{ {./eps/pdf/} }
\tikzset{every loop/.style={}}
\usepackage[pagewise]{lineno} 
\usepackage{cleveref}

\newtheorem{prop}{Proposition}

\newtheorem{defi}{Definition}
\newtheorem{con}{Conjecture}
\newtheorem{theo}{Theorem}
\theoremstyle{definition}

\theoremstyle{remark}
\newtheorem{rem}{Remark}

\crefdefaultlabelformat{#2{\scshape #1}#3}
\crefname{subfigure}{fig.}{figs.}
\Crefname{subfigure}{Figure}{Figure}

\title[SYMMETRY BREAKING IN A GLOBALLY COUPLED MAP OF
FOUR SITES
]{SYMMETRY BREAKING IN A GLOBALLY COUPLED MAP OF
FOUR SITES
}
\author[F. M. S\'elley]{Fanni M. S\'elley$^1$}
\address{$^1$Department of Stochastics, Institute of Mathematics \\
Budapest University of Technology and Economics \\
Egry J\'{o}zsef u. 1, H-1111 Budapest, Hungary
and
MTA-BME Stochastics Research Group \\
Budapest University of Technology and Economics \\
Egry J\'ozsef u. 1, H-1111 Budapest, Hungary.
}
\date{\today}

\begin{document}

\begin{abstract}
A system of four globally coupled doubling maps is studied in this paper.
It is known that such systems have a unique absolutely continuous
invariant measure (acim) for weak interaction, but the case of stronger coupling is still unexplored. As in the case of three coupled sites \cite{selley2015meanfield}, we prove the existence of a critical value
of the coupling parameter at which multiple
acims appear. Our proof has several new ingredients in comparison to the one presented in
\cite{selley2015meanfield}. We strongly exploit the symmetries of the dynamics in the course of the argument.
This simplifies the computations considerably, and gives us a precise description of the
geometry and symmetry properties of the arising asymmetric invariant sets. Some new
phenomena are observed which are not present in the case of three sites. In particular, the asymmetric
invariant sets arise in areas of the phase space which are transient for weaker coupling and a
nontrivial symmetric invariant set emerges, shaped by an underlying centrally symmetric
Lorenz map. We state some conjectures on further invariant sets, indicating that unlike the case of three sites, ergodicity breaks down in many steps, and not all of them are
accompanied by symmetry breaking.

\end{abstract}

\maketitle

\let\thefootnote\relax\footnotetext{\emph{AMS subject classification.} 37A25; 37E10; 37G99.}
\let\thefootnote\relax\footnotetext{\emph{Key words and phrases.} coupled map systems, ergodicity breaking, asymmetric invariant sets, piecewise
affine maps.}

\section{Introduction}

Coupled map systems are simple models of a finite or infinite network of interacting units generally referred to as sites. The dynamics is given by the composition of the (typically chaotic) individual dynamics and a coupling map representing the characteristics of the interaction. The coupling map usually includes a parameter $0 \leq \varepsilon \leq 1$, representing the strength of interaction. 

The analysis of coupled map systems is quite a challenging task. The main interest undoubtedly lies in the emergence of bifurcations: as the strength of interaction varies, the main features of coupled maps can change dramatically. To give a brief overview of the existing literature, we only cite papers connected closely to our work. For more complete lists, see the references in \cite{koiller2010coupled}, \cite{fernandez2014breaking}, \cite{keller2006uniqueness} and the collection \cite{chazottes2005dynamics}. Most results concern the case of coupling strength close to zero. In this case, the behavior of the coupled system is, in general, similar to that of the uncoupled system. For example, the existence of a unique SRB measure, possibly with strong chaotic properties can be proved  (\cite{jarvenpaa1997srb},\cite{jiang1998equilibrium}, \cite{keller2006uniqueness}). Further results are related to the emergence of contracting directions for values of $\varepsilon$ close to 1, resulting in the absence of an absolutely continuous invariant measure (\cite{boldrighini1995ising}, \cite{boldrighini2001ising}, \cite{just1995globally},\cite{koiller2010coupled}). The results for higher values of the coupling strength can be usually thought of as synchronization in some simple sense. Examples include the individual systems behaving asymptotically identically (see the case of two coupled maps for $\varepsilon > 1/2$ in \cite{selley2015meanfield}) or acquiring some more complicated, yet fixed formation (three systems acting interdependently on the circle are shown to be able to acquire evenly placed positions asymptotically for $\varepsilon > 1/2$ in \cite{koiller2010coupled}).  

Furthermore, numerical simulations suggest that more complicated phenomena are possible. A particularly interesting one is the emergence of \emph{multiple absolutely continuous invariant measures} (acims) in the case when the system is still fully expanding. Such bifurcations can be interpreted as a deterministic analogue of the phase transitions of Ising models in statistical physics (\cite{bunimovich1988spacetime}, \cite{gielis2000coupled},\cite{miller1993macroscopic}). The coupling parameter in this case should be relatively high, so that the previously mentioned perturbative results do not apply. Hence these results might be thought of as complex synchronization phenomena. For example, see the case of three coupled maps in \cite{selley2015meanfield}, where it was shown that ergodic components could be interpreted in terms of the relative positions of the sites.

In this paper we are going to study a system specifically constructed to demonstrate such phenomena, often termed ergodicity breaking. The model, introduced by Koiller and Young in \cite{koiller2010coupled}, is a globally coupled system of $N$ identical circle maps.  By standard results in the literature \cite{keller2006uniqueness}, the map has a unique mixing acim for $\varepsilon$ values close to zero. Fernandez \cite{fernandez2014breaking} indicated, by numerically calculating certain order parameters, that multiple acims emerge when the value of the coupling parameter is increased. In other words, ergodicity is broken at some critical value of $\varepsilon$. 

Precise results exist for a small number of sites. Fernandez showed in \cite{fernandez2014breaking} that ergodicity breaking does not occur for $N=2$. Bálint and the present author pointed out later that even though the unique acim is ergodic in the expanding regime, it ceases to be mixing if $1-\frac{\sqrt{2}}{2} \leq \varepsilon$ \cite{selley2015meanfield}. In the case of $N=3$, it was shown by Fernandez and independently by Bálint and the present author, that due to the appearance of asymmetric invariant sets, ergodicity breaking occurs at the value $\frac{4-\sqrt{10}}{2} \approx 0.42$ of the coupling parameter.

However, for $N > 3$, only numerical simulations suggest ergodicity breaking and no rigorous results have been obtained yet. The existing results for $N=2$ and $N=3$ were acquired with the help of elementary, yet careful geometric considerations, requiring a detailed understanding of the dynamics. This paper considers \emph{the case of $N=4$} in this spirit. The results for $N=2$ were straightforward consequences of well-known facts about centrally symmetric Lorenz maps. For $N=3$, it was not possible to use known results, and a precise geometric understanding of the action of a certain 2-dimensional map was necessary. The case is similar with $N=4$, but now the geometry is considerably more complicated, since the understanding of a 3-dimensional map is required. B. Fernandez revealed in private communications, that he identified two values of the coupling parameter, $\varepsilon_a \approx 0.39$ and $\varepsilon_b \approx 0.43$, at which numerical simulations indicate the appearance of multiple asymmetric invariant sets. He also devised an algorithm, which by fixing $\varepsilon > \varepsilon_a$, generates an asymmetric invariant set of the system with coupling strength $\varepsilon$. This algorithm could act as a base for a \emph{computer assisted proof} for ergodicity breaking. Nevertheless, it leaves us with some questions which an analytic proof could answer. For example, the set generated by the algorithm is difficult to interpret geometrically, and the algorithm does not provide any rigorous results about the symmetry properties of the set. The main goal of this paper is to present an \emph{analytic proof of ergodicity breaking}. Although the result is similar to the one obtained in the $N=3$ case in \cite{selley2015meanfield}, we now worked out a more systematic proof which avoids many of the redundancies the previous proof exhibited. In particular, we exploit the symmetries of the system and use simple facts of linear optimization to decrease the amount of calculations as much as possible in this way. This methodological simplification is essential given that the $N=4$ case is far more complex. 

The paper is organized as follows: in Section \ref{SecDyn}, we familiarize ourselves with the dynamics and apply the change of coordinates introduced in \cite{selley2015meanfield}, to reduce the analysis of our original 4-dimensional system to the study of a piecewise affine map of the 3-dimensional torus. We describe the system defined by this map in detail and explore its symmetry group. Section \ref{SecMain} contains our theorem stating that there exists a critical value $\varepsilon^*$ such that the system is not ergodic for $\varepsilon^* \leq \varepsilon$.  In Section \ref{SecInv} we give the technical details needed to prove our theorem and state some conjectures explicitly. In particular, we describe an asymmetric set $\mathcal{A}$ that is invariant for $\varepsilon^* \leq \varepsilon$ and does not break a special abelian subgroup of the symmetry group. From simulations it seems clear that this set is in a part of the phase space which is transient for $1-\frac{\sqrt{2}}{2} \leq \varepsilon < \varepsilon^*$, so $\mathcal{A}$ does not appear as a decomposition of a symmetric invariant set, in contrast to the $N=3$ case. We further describe a nontrivial symmetric invariant set $\mathcal{S}$. Our conjectures in this section concern the appearance of further invariant sets. We conjecture that another asymmetric invariant set appears at some value $\varepsilon^{**} > \varepsilon^*$. Furthermore, the set $\mathcal{S}$ breaks up countably many times into new symmetric invariant sets as $\varepsilon$ goes to $\frac{1}{2}$. Based solely on the $N=3$ case it seemed that ergodicity breaks down in a single step. Now our results and observations indicate that the situation is more complex for higher values of $N$. Although the proofs are straightforward, they are lengthy, and thus left to the appendices. Appendix A contains important reference for Appendices B and C, which contain the proofs of our statements regarding the sets $\mathcal{A}$ and $\mathcal{S}$, respectively. 


\section{Definition of the dynamics} \label{SecDyn}
\begin{figure}[h!]
 \centering
 \begin{tikzpicture}[scale=2]
       \draw[->] (-2.5,0) -- (2.5,0) node[above] {$u$};
       \draw[->] (0,-1) -- (0,1) node[right] {\hspace{0.1cm}$g(u)$};
       \draw[very thick,red] (-0.7,-0.7) -- (0.7,0.7);
       \draw[very thick,red] (0.7,-0.7) -- (2.1,0.7);
       \draw[very thick,red] (-0.7,0.7) -- (-2.1,-0.7);
       \draw[red, fill] (0.7,0) circle (0.03cm);
       \draw[red, fill] (-0.7,0) circle (0.03cm);
       \draw[red, fill] (2.1,0) circle (0.03cm);
       \draw[red, fill] (-2.1,0) circle (0.03cm);
       \foreach \x/\xtext in {0.7/\frac{1}{2},1.4/1, 2.1/\frac{3}{2},-0.7/-\frac{1}{2},-1.4/-1, -2.1/-\frac{3}{2}}
           \draw[shift={(\x,0)}] (0pt,2pt) -- (0pt,-2pt) node[below] {$\xtext$};
       \foreach \y/\ytext in {0.7/\frac{1}{2},-0.7/-\frac{1}{2}}
             \draw[shift={(0,\y)}] (2pt,0pt) -- (-2pt,0pt) node[left] {$\ytext$};
     \end{tikzpicture}
     \caption{The function $g$.} \label{g}
     \end{figure}
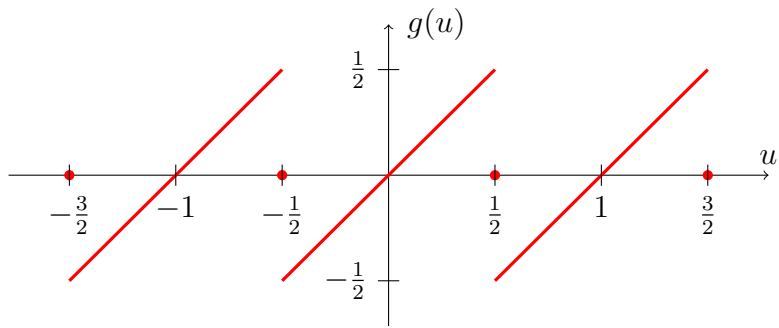

We are going to consider the system of globally coupled doubling maps as defined by Fernandez in \cite{fernandez2014breaking}. Let $N>0$ and $\mathbb{T}^N=(\mathbb{R}\backslash \mathbb{Z})^N$ be the $N$-dimensional torus. We are going to use the representation $\mathbb{T} \equiv [0,1]$, hence $\mathbb{T}^N \equiv [0,1]^N$. We define $F_{\varepsilon,N}: \mathbb{T}^N \to \mathbb{T}^N$ as
\begin{equation}
(F_{\varepsilon,N}(x))_i=2\left(x_i+\frac{\varepsilon}{N} \sum_{j=1}^N g(x_j-x_i)\right) \quad \mod 1, \hspace{0.3cm} x=(x_s)_{s=1}^N \in \mathbb{T}^N, i=1\dots N,
\end{equation} 
where the function $g$ responsible for the features of the interaction is defined as the lift of
\[
\hat{g}(u)=\begin{cases}
0 & \text{ if } u = \pm \frac{1}{2} , \\
u& \text{ if } u \in \left(-\frac{1}{2}, \frac{1}{2} \right)
\end{cases}
\]
to $\mathbb{R}$, see Figure \ref{g} for the graph of this map.
 
The map $F_{\varepsilon,N}$ can be regarded as the composition of a mean-field type coupling map
\[
(\Phi_{\varepsilon,N}(x))_i=x_i+\frac{\varepsilon}{N} \sum_{j=1}^N g(x_j-x_i) \quad \mod 1, \hspace{0.3cm} x=(x_s)_{s=1}^N \in \mathbb{T}^N, i=1\dots N,
\] 
and the individual dynamics $T(x)=2x \mod 1$, a map well known for its strong chaotic properties. 

This dynamical system can be thought of as a model of $N$ interacting particles on the circle, changing position according to both the individual dynamics $T$ and the coupling map $\Phi_{\varepsilon,N}$, representing the effect of the other particles. The parameter $\varepsilon \in [0,1/2)$ is called the coupling parameter, and it controls the strength of the interaction between particles. In particular, $\varepsilon=0$ means no interaction, the position of the particles evolve independently according to $T$. We only consider values lower than 1/2 for the coupling parameter, since this is the expanding regime of the map $F_{\varepsilon,N}$ (more precisely, this is the range of the coupling parameter where all of the eigenvalues of the Jacobian of $F_{\varepsilon,N}$ are greater than 1.) 

It follows from the abstract framework of \cite{keller2006uniqueness} that if $\varepsilon$ is sufficiently small, this system has a unique absolutely continuous invariant measure. The main goal of our upcoming analysis is to show that in the case of $N=4$, if $\varepsilon$ is greater than some critical value (which is smaller than $1/2$), this ceases to be true. 

So from now on we restrict our attention to the special case of four sites. The dynamics of the 4-dimensional system, $F_{\varepsilon,4}: \mathbb{T}^4 \to \mathbb{T}^4$, takes the form
\begin{align*}
F_{\varepsilon,4}(x_1,x_2,x_3,x_4)=\bigg(&2x_1+\frac{\varepsilon}{2}\left(g(x_2-x_1)+g(x_3-x_1)+g(x_4-x_1)\right), \\
&2x_2+\frac{\varepsilon}{2}\left(g(x_1-x_2)+g(x_3-x_2)+g(x_4-x_2)\right), \\
&2x_3+\frac{\varepsilon}{2}\left(g(x_1-x_3)+g(x_2-x_3)+g(x_4-x_3)\right), \\
&2x_4+\frac{\varepsilon}{2}\left(g(x_1-x_4)+g(x_2-x_4)+g(x_3-x_4)\right)\bigg) \qquad \mod 1.
\end{align*}

in the coordinates $x_1,x_2,x_3$ and $x_4$.

As in \cite{selley2015meanfield}, we define the following new coordinates:
\begin{table}[h!!]
\centering
\begin{tabular}{lll}
$s$&$=x_1+x_2+x_3+x_4$ & $\mod 1,$\\
$p$&$=x_1-x_2$ & $\mod 1,$\\
$q$&$=x_2-x_3$ & $\mod 1,$\\ 
$r$&$=x_3-x_4$ & $\mod 1.$
\end{tabular}
\end{table}

Let us consider
\[
T \times G_{\varepsilon,3}: \mathbb{T} \times \mathbb{T}^3 \to \mathbb{T} \times \mathbb{T}^3,
\]
where $G_{\varepsilon,3}$ is the following map of $\mathbb{T}^3$:
\begin{align*}
G_{\varepsilon,3}(p,q,r)=\bigg(&2p+\frac{\varepsilon}{2}\left(-2g(p)+g(q)+g(q+r)-g(p+q)-g(p+q+r)\right), \\
&2q+\frac{\varepsilon}{2}\left(-2g(q)+g(p)+g(r)-g(q+r)-g(p+q)\right), \\
&2r+\frac{\varepsilon}{2}\left(-2g(r)+g(q)+g(p+q)-g(q+r)-g(p+q+r)\right)\bigg) \qquad \mod 1.
\end{align*}
It is important to note that the system with law $T \times G_{\varepsilon,3}$ is not conjugate to the system with law $F_{\varepsilon,4}$, but a factor of it, since the points 
\[
\left(x_1+\frac{i}{4},x_2+\frac{i}{4},x_3+\frac{i}{4},x_4+\frac{i}{4}\right), \quad i=1\dots 4 
\]
share the same $s,p,q,r$-coordinates. However, detection of ergodicity breaking in a factor implies ergodicity breaking in the original system. Furthermore, ergodicity breaking in the system
\[
G_{\varepsilon,3} : \mathbb{T}^3 \to \mathbb{T}^3
\]
implies the lack of ergodicity of the factor, so we are going to continue with the analysis of this system.

We are going to represent the domain of $G_{\varepsilon,3}$, that is $\mathbb{T}^3$, as the unit cube in $\mathbb{R}^3$ with opposite faces identified. The map $G_{\varepsilon,3}$ is piecewise affine, and the singularities arise from the singularities of the function $g$, giving polyhedral domains of continuity. For a complete description of the singularities, the continuity domains and the precise action of $G_{\varepsilon,3}$ on them, see Appendix A. 

In the next subsection we discuss the symmetries of this map.

\subsection{Symmetries of the map $G_{\varepsilon,3}$} \label{SecSym}

A symmetry of a map $F$ is a linear transformation $S$ such that 
\[
S\circ F=F \circ S.
\]
The symmetries of the map $F_{\varepsilon,4}$ arise from two sources:
\begin{itemize}
\item the inversion symmetry of $g$ and that of the doubling map (namely, $g(1-u)=1-g(u)$ and $T(1-u)=1-T(u)$) imply the inversion symmetry of $F_{\varepsilon,4}$. More precisely, 
\[
F_{\varepsilon,4} \circ I = I \circ F_{\varepsilon,4},
\]
where
\[
I: (x_1,x_2,x_3,x_4) \to (1-x_1,1-x_2,1-x_3,1-x_4).
\]
\item every permutation of $x_1,x_2,x_3,x_4$ is a symmetry of $F_{\varepsilon,4}$:
\[
F_{\varepsilon,4} \circ \pi = \pi \circ F_{\varepsilon,4},
\]
where $\pi$ is an element of the fourth order symmetric group (the group of all permutations of four elements).
\end{itemize}

We note that the symmetry group of the map $F_{\varepsilon,4}$ is generated by the inversion $I$ and a generator of the fourth order symmetric group. An example for the latter is
\begin{align*}
 \pi_1:& \hspace{0.2cm}(x_1,x_2,x_3,x_4) \mapsto(x_2,x_1,x_3,x_4), \\
 \pi_2:& \hspace{0.2cm} (x_1,x_2,x_3,x_4) \mapsto(x_3,x_2,x_1,x_4), \\
 \pi_3:& \hspace{0.2cm} (x_1,x_2,x_3,x_4) \mapsto(x_4,x_2,x_3,x_1), \\
 \pi_4:& \hspace{0.2cm} (x_1,x_2,x_3,x_4) \mapsto(x_1,x_3,x_2,x_4), \\
 \pi_5:& \hspace{0.2cm} (x_1,x_2,x_3,x_4) \mapsto(x_1,x_4,x_3,x_2), \\
 \pi_6:& \hspace{0.2cm} (x_1,x_2,x_3,x_4) \mapsto(x_1,x_2,x_4,x_3).
\end{align*}

Note that this is not a minimal generator in the sense that for example $\pi_1,\pi_2$ and $\pi_3$ already generate the symmetric group. However, we are working with a mean-field model, and no coordinate can have a special role. We constructed a generator with as few elements as possible such that every coordinate has the same role. 

These symmetries induce the generators of the symmetry group of $G_{\varepsilon,3}$, which we shall denote by $S_G$. The inversion of $(x_1,x_2,x_3,x_4)$ induces the inversion of $(p,q,r)$:
\[
S_0: \quad (p,q,r) \mapsto (1-p,1-q,1-r)
\]
The permutations induce the following symmetries:

\begin{table}
\centering
\begin{tabular}{lll}
$S_1:$ &  $\begin{bmatrix}p \\
q \\
r \end{bmatrix} \mapsto \begin{bmatrix}-1 & 0 & 0 \\
1 & 1 & 0 \\
0 & 0 & 1 \end{bmatrix} \begin{bmatrix}p \\
q \\
r \end{bmatrix}$ & $\mod 1$, \\
\rule{0pt}{6ex}
$S_2:$ &  $\begin{bmatrix}p \\
q \\
r \end{bmatrix} \mapsto \begin{bmatrix}0 & -1 & 0 \\
-1 & 0 & 0 \\
1 & 1 & 1 \end{bmatrix} \begin{bmatrix}p \\
q \\
r \end{bmatrix}$ & $\mod 1$, \\
\rule{0pt}{6ex}
$S_3:$&  $\begin{bmatrix}p \\
q \\
r \end{bmatrix} \mapsto \begin{bmatrix}0 & -1 & -1 \\
0 & 1 & 0 \\
-1 & -1 & 0 \end{bmatrix} \begin{bmatrix}p \\
q \\
r \end{bmatrix}$ & $\mod 1$, \\
\rule{0pt}{6ex}
$S_4:$&  $\begin{bmatrix}p \\
q \\
r \end{bmatrix} \mapsto \begin{bmatrix}1 & 1 & 0 \\
0 & -1 & 0 \\
0 & 1 & 1 \end{bmatrix} \begin{bmatrix}p \\
q \\
r \end{bmatrix}$ & $\mod 1$, \\
\rule{0pt}{6ex}
$S_5:$&  $\begin{bmatrix}p \\
q \\
r \end{bmatrix} \mapsto \begin{bmatrix}1 & 1 & 1 \\
0 & 0 & -1 \\
0 & -1 & 0 \end{bmatrix} \begin{bmatrix}p \\
q \\
r \end{bmatrix}$ & $\mod 1$, \\
\rule{0pt}{6ex}
$S_6:$&  $\begin{bmatrix}p \\
q \\
r \end{bmatrix} \mapsto \begin{bmatrix}1 & 0 & 0 \\
0 & 1 & 1 \\
0 & 0 & -1 \end{bmatrix} \begin{bmatrix}p \\
q \\
r \end{bmatrix}$ & $\mod 1$.
\end{tabular}
\end{table}

Note that
\begin{align*}
\pi_4&=\pi_2\pi_1\pi_2 \quad \Rightarrow \quad S_4=S_2S_1S_2, \\
\pi_5&=\pi_3\pi_1\pi_3 \quad \Rightarrow \quad S_5=S_3S_1S_3, \\
\pi_6&=\pi_3\pi_2\pi_3 \quad \Rightarrow \quad S_6=S_3S_2S_3. 
\end{align*}
We further note two facts: since $I$ commutes with every symmetry of $F_{\varepsilon,4}$, $S_0$ commutes with every symmetry of $G_{\varepsilon,3}$. Also, each generator symmetry is a $\mathbb{Z}_2$-symmetry, applying it twice yields identity.

\section{Main result} \label{SecMain}

We first define the notions of symmetry and invariance we are going to use. 
\begin{defi}
A set $\mathcal{B} \subset \mathbb{T}^3$ is symmetric with respect to $S \in S_G$ if $\mathcal{B}=S\mathcal{B}$, and asymmetric with respect to $S \in S_G$ if $\mathcal{B}$ and its symmetric image $S\mathcal{B}$ are disjoint.

A set $\mathcal{B}$ is symmetric if it is symmetric with respect to every element of $S_G$, and asymmetric if there exists a symmetry in $S_G$ for which $\mathcal{B}$ is asymmetric.
\end{defi}

To keep terminology brief, we are going to say that $\mathcal{B}$ \emph{breaks} $S$, when $\mathcal{B}$ is asymmetric with respect to some $S$.

\begin{defi}
A set $\mathcal{B} \subset \mathbb{T}^3$ is (forward) invariant if $G_{\varepsilon,3}(\mathcal{B}) \subseteq \mathcal{B}$.  
\end{defi}

Notice that if $\mathcal{B}$ is invariant under the dynamics, $S\mathcal{B}$ is also invariant, since
\[
G_{\varepsilon,3}(\mathcal{B}) \subseteq \mathcal{B} \Rightarrow G_{\varepsilon,3}(S\mathcal{B})= SG_{\varepsilon,3}(\mathcal{B}) \subseteq S\mathcal{B}.
\]

Suppose this set $\mathcal{B}$ is a (not necessarily connected) polyhedral region of $\mathbb{T}^3$ (this will be the relevant case for us) and asymmetric with respect to some $S \in S_G$. Of course in this case the symmetric images of $\mathcal{B}$ are also polyhedral domains and disjoint from $\mathcal{B}$. Since the map is completely expanding (as $\varepsilon < 1/2$), we can use Theorem 1.7 in \cite{thomine2010spectral} for these sets independently, to obtain that on each of these sets an absolutely continuous invariant measure is supported. In conclusion, the existence of such an asymmetric invariant set means that multiple acims exist. Or from another point of view, the acim of maximal support cannot be ergodic.

In our theorem we state that if the value of the coupling parameter is sufficiently large, a set like this exists.

\begin{theo}
There exists an $\varepsilon^* < 1/2$, such that for $\varepsilon^* \leq \varepsilon < 1/2$ the system $G_{\varepsilon,3}: \mathbb{T}^3 \to \mathbb{T}^3$ admits an asymmetric invariant set $\mathcal{A}$. 
\end{theo}

The proof of Theorem 1 follows from Proposition 1 of Section \ref{SecAsym}, which includes a precise value of such an $\varepsilon^*$, conjectured to be the smallest.  Proposition 1 also includes the definition of the set $\mathcal{A}$ and a description of its symmetry properties. 

At the end of Section \ref{SecAsym}, we state a conjecture on the existence of $\varepsilon^{**}$ larger then $\varepsilon^{*}$, such that another asymmetric invariant set exists provided that $\varepsilon \geq \varepsilon^{**}$. We describe this set and a probable value of $\varepsilon^{**}$ in Conjecture 1.

This would indicate that the number of acims increases not only at the value of $\varepsilon^*$, but also at another larger value $\varepsilon^{**}$. In Section \ref{SecSymG} we formulate a conjecture that countably many further critical values of the coupling parameter exist corresponding to the appearance of new invariant sets, but the new sets are most likely symmetric. These sets are all contained in a special symmetric invariant set, described in Proposition 2. We give a conjecture for the countably many critical values of the coupling parameter in Conjecture 2.   

\section{Invariant sets} \label{SecInv}

\subsection{A centrally symmetric Lorenz map} \label{SecL}
In this subsection we describe an interval map, which is going to play a central role in the definitions of the invariant sets to be described in this section. Let us define $L_{\varepsilon}: [\varepsilon/2, 1-\varepsilon/2] \to [\varepsilon/2, 1-\varepsilon/2]$ as
\begin{equation} \label{h}
  L_{\varepsilon}(v)=
  \begin{cases}
     2(1-\varepsilon)v+\frac{\varepsilon}{2} & \text{if }  \frac{\varepsilon}{2} < v < \frac{1}{2}, \\
     2(1-\varepsilon)v+\frac{3\varepsilon}{2}-1       & \text{if } \frac{1}{2} < v < 1-\frac{\varepsilon}{2}.
    \end{cases}
  \end{equation}
The graph of this map is plotted on Figure \ref{h}.

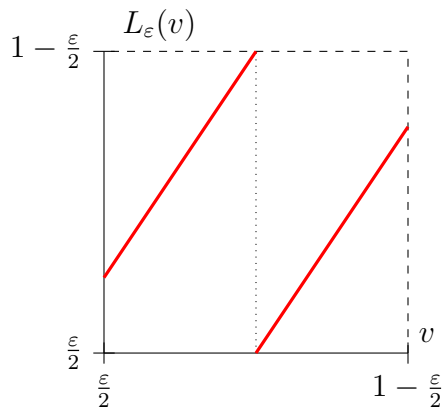
\begin{figure}[h!]
 \centering
 \begin{tikzpicture}[scale=2]
 \draw (0,0) -- (2,0) node[above right] {$v$};
        \draw (0,0) -- (0,2) node[above right] {\hspace{0.1cm}$L_{\varepsilon}(v)$};
        \draw[dashed] (2,0) -- (2,2) -- (0,2);
        \draw[dotted] (1,0) -- (1,2);
        \draw[very thick,red] (0,0.5) -- (1,2);
        \draw[very thick,red] (1,0) -- (2,1.5);
        \foreach \x/\xtext in {2/1-\frac{\varepsilon}{2}, 0/\frac{\varepsilon}{2}}
            \draw[shift={(\x,0)}] (0pt,2pt) -- (0pt,-2pt) node[below] {$\xtext$};
        \foreach \y/\ytext in {2/1-\frac{\varepsilon}{2},0/\frac{\varepsilon}{2}}
              \draw[shift={(0,\y)}] (2pt,0pt) -- (-2pt,0pt) node[left] {$\ytext$};
 \end{tikzpicture}
 \caption{The graph of $L_{\varepsilon}$.} \label{hl}
 \end{figure}  

This is a well known centrally symmetric Lorenz map. We are going to state some known facts which will prove useful for us, following the classical work of Parry \cite{parry1979lorenz}.

For every $0 \leq \varepsilon < 1/2$, the map has an ergodic invariant measure, supported on a finite union of intervals. If $\varepsilon < 1-\frac{\sqrt{2}}{2}$, the invariant measure is supported on one interval, and it is mixing. For every integer $n$, when 
\[
1-\frac{\sqrt[2^n]{2}}{2} \leq \varepsilon <  1-\frac{\sqrt[2^{n+1}]{2}}{2},
\] 
the supporting intervals of the invariant measure can be grouped in $2^n$ mixing components. Each mixing component is a union of intervals, restricted to which $L_{\varepsilon}^{2^n}$ has a mixing invariant measure. 

For example, take $n=1$, yielding
\[
1-\frac{\sqrt{2}}{2} \leq \varepsilon <  1-\frac{\sqrt[4]{2}}{2}.
\]
Two mixing components exist for these values of $\varepsilon$:
 \begin{equation} \label{mixing}
 (\varepsilon/2, L_{\varepsilon}^2(1-\varepsilon/2)) \cup ( L_{\varepsilon}^2(\varepsilon/2), 1-\varepsilon/2) \quad \text{ and } \quad  (L_{\varepsilon}(\varepsilon/2), L_{\varepsilon}(1-\varepsilon/2)).
 \end{equation}
These two components are mapped to one another by $L_{\varepsilon}$ and the restriction of $L^2_{\varepsilon}$ to each of them has a mixing invariant measure. For further details, see Parry \cite{parry1979lorenz}. 

Lastly, we state that the map $L_{\varepsilon}$ has a 2-periodic orbit $p^* \leftrightarrow 1-p^*$, where
\begin{equation} \label{pstar}
p^*=\frac{(1-\varepsilon)\varepsilon+3\varepsilon/2-1}{1-4(1-\varepsilon)^2}=\frac{\varepsilon-2}{4\varepsilon-6}.
\end{equation}

 \subsection{Asymmetric invariant sets} \label{SecAsym}
 
 \begin{figure}[h!]
    \centering
    \includegraphics[scale=0.6]{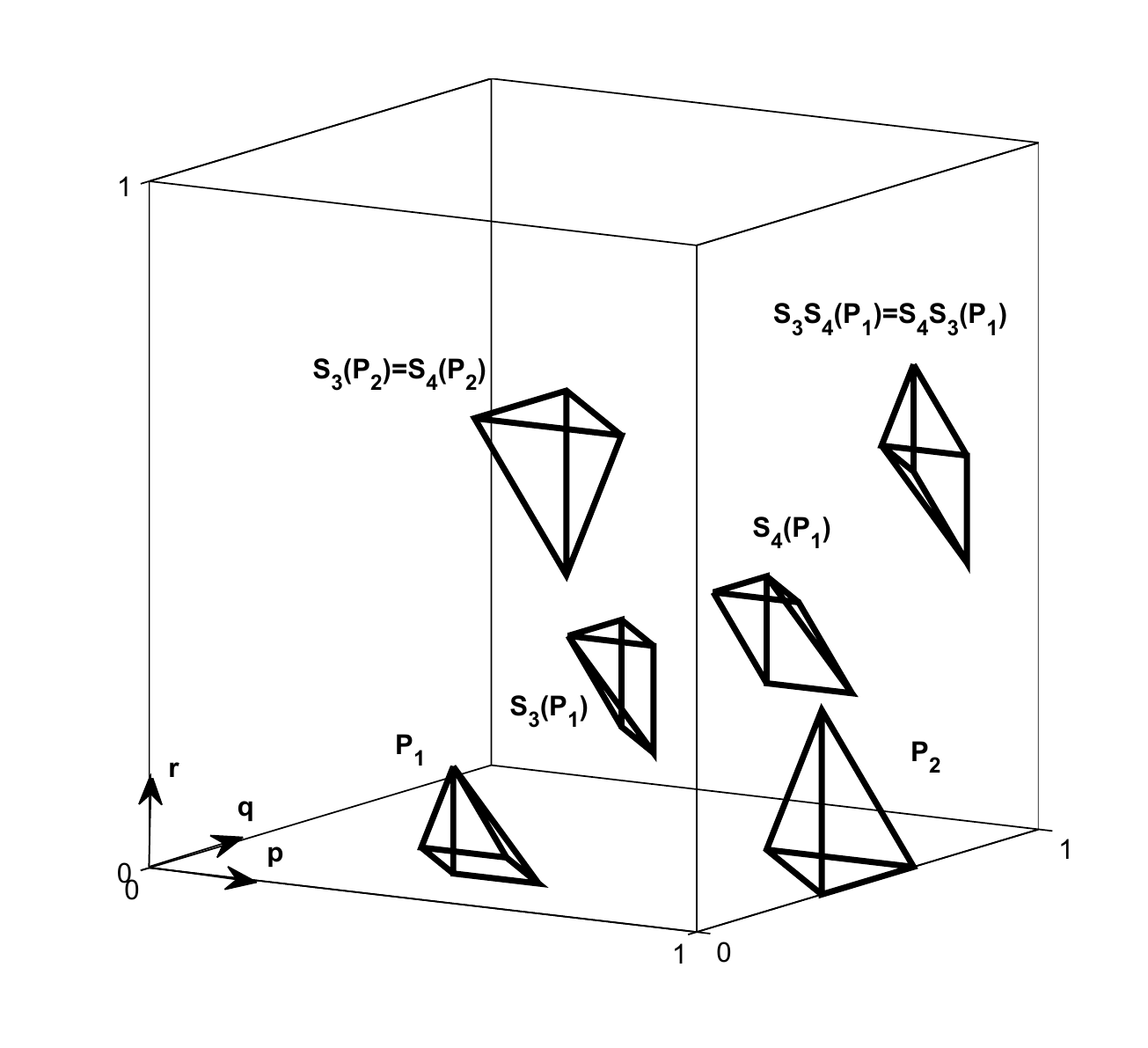}
    \caption{The asymmetric set $\mathcal{A}$.} \label{FigA}
    \end{figure}
    
 In this section we define a value $\varepsilon^*$ of the coupling parameter, and give an asymmetric invariant set of the map $G_{\varepsilon,3}$ when $\varepsilon^* \leq \varepsilon < 1/2$. Before describing the set explicitly, we are going to give a few words about our intuition leading to this particular set.
 
 It is a somewhat natural thought that the set should be the union of polyhedra with faces parallel to singularities. Using that on certain invariant circles the restricted map is exactly $L_{\varepsilon}$, one can explicitly calculate 12 period two points of the map $G_{\varepsilon,3}$. We choose the faces of the polyhedra so that the period two points of $G_{\varepsilon,3}$ are some of their vertices. This is very similar to the construction of asymmetric invariant sets in the $N=3$ case. This means that the constants defining the planes which give the faces of the polyhedra include the period two point $p^*$ of $L_{\varepsilon}$, defined in terms of $\varepsilon$ by equation \eqref{pstar}. Actually, the exact position of the planes will be given either by the constant $p^*, 1-p^*$ mod 1 or $\varepsilon/2, 1-\varepsilon/2$ mod 1 (these latter constants arise from the first iteration image of the continuity domains, which give a natural constraint on a forward invariant set). In turns out that a clever such choice ensures invariance for large enough $\varepsilon$. To ensure asymmetry, we define the set to be invariant under a proper subgroup of the symmetry group. 
 
 So let us define two polyhedra, $P_1$ and $P_2$. We note here that we are going to define polyhedra with a set of linear inequalities which give a minimal representation in $\mathbb{R}^3$. Let $P_1$ and $P_2$ be defined in the following way: 
 
 \begin{center}
     \begin{tabular}{|p{2cm}|p{5cm}|p{5cm}|}
       \hline 
       & $\mathbf{P_1}$ & $\mathbf{P_2}$  \\ \hline
       $\mathbf{p}$&&$p < 1$ \\ \hline 
       $\mathbf{q}$&$q > \varepsilon/2$&\\ \hline 
       $\mathbf{r}$&$r > 0$&$r > 0$ \\ \hline
       $\mathbf{p+q}$&$p+q > 1-p^*$&$p+q > 1+p^*$ \\ \hline
       $\mathbf{q+r}$&$q+r < p^*$&$q+r < 1-p^*$ \\ \hline  
       $\mathbf{p+q+r}$&$p+q+r < 1-\varepsilon/2$& \\ \hline 
     \end{tabular}
 \end{center}

 \begin{prop} \label{propestar1}
 The set
 $$\mathcal{A}=P_1 \cup P_2 \cup S_3(P_1) \cup S_4(P_1) \cup S_3S_4(P_1) \cup S_3(P_2)$$
 is invariant with respect to $G_{\varepsilon,3}$ if and only if
 $$0.397 \approx \varepsilon^* =\frac{1}{6}\left(7-4\sqrt[3]{\frac{1}{43-3\sqrt{177}}}-\sqrt[3]{\frac{43-3\sqrt{177}}{2}} \right) \leq  \varepsilon $$
 and it is symmetric with respect to $S_3$ and $S_4$, but asymmetric with respect to
 $$S_0, \text{ } S_1, \text{ } S_2, \text{ } S_5, \text{ and } S_6.$$
 \end{prop}


\begin{rem}
The precise value of $\varepsilon^*$ may seem complicated. It is actually the unique real solution of
\[
p^*=(1-\varepsilon)^2, \quad \text{or more explicitly} \quad  4\varepsilon^3-14\varepsilon^2+15\varepsilon-4=0.
\]
The importance of this equation is explained in Appendix B, see equations \eqref{epsstar1}-\eqref{epsstar2} and the argument leading up to them.  
\end{rem}
 
 \begin{rem}
 \begin{figure}[h!!]
  \centering
  \begin{tikzpicture}[>=stealth',shorten >=1pt,auto,node distance=4cm,
                  very thick,main node/.style={circle,draw,font=\L_{\varepsilon}arge\bfseries},scale=0.95]
  
    \node[draw=none,fill=none] (00) { $\mathcal{A}$};
    \node[draw=none,fill=none] (1) at (-4,-2) {$S_1(\mathcal{A})$};
    \node[draw=none,fill=none] (0) at (4,-2) {$S_0(\mathcal{A})$};
    \node[draw=none,fill=none] (2) [below of=1] {$S_2(\mathcal{A})$};
    \node[draw=none,fill=none] (6) [below of=0] {$S_6(\mathcal{A})$};
    \node[draw=none,fill=none] (5) at (0,-8) {$S_5(\mathcal{A})$};
  
    \path
      (00) edge [black, out=160,in=110,looseness=8] node {} (00)
      (00) edge [yellow, out=70,in=20,looseness=8] node {} (00)
      (00) edge [blue] node {} (1)
      (00) edge [red] node {} (2)
      (00) edge [cyan] node {} (5)
      (00) edge [magenta] node {} (6)
      (00) edge [green] node {} (0)
      
      (1) edge [red, out=160,in=100,looseness=3] node {} (1)
      (1) edge [cyan, out=260,in=200,looseness=3] node {} (1)
      (1) edge [yellow] node {} (2)
      (1) edge [black] node {} (5)
      (1) edge [magenta] node {} (0)
      (1) edge [green] node {} (6)
      
      (2) edge [magenta, out=280,in=230,looseness=6] node {} (2)
      (2) edge [blue, out=190,in=140,looseness=4] node {} (2)
      (2) edge [black] node {} (6)
      (2) edge [green] node {} (5)
      (2) edge [cyan] node {} (0)
      
      (5) edge [magenta, out=240,in=190,looseness=3] node {} (5)
      (5) edge [blue, out=350,in=300,looseness=3] node {} (5)
      (5) edge [yellow] node {} (6)
      (5) edge [red] node {} (0)
      
      (6) edge [red, out=260,in=310,looseness=6] node {} (6)
      (6) edge [cyan,  out=350,in=40,looseness=4] node {} (6)
      (6) edge [blue] node {} (0)
      
      (0) edge [black, out=20,in=80,looseness=3] node {} (0)
      (0) edge [yellow, out=280,in=340,looseness=3] node {} (0);

  \end{tikzpicture}
  \caption{Six asymmetric invariant sets and the symmetries connecting them. Edge colors indicate the following symmetries: green -- $S_0$, blue -- $S_1$, red -- $S_2$, black -- $S_3$, yellow -- $S_4$, cyan -- $S_5$, magenta -- $S_6$.} \label{FigAsym}
  \end{figure}
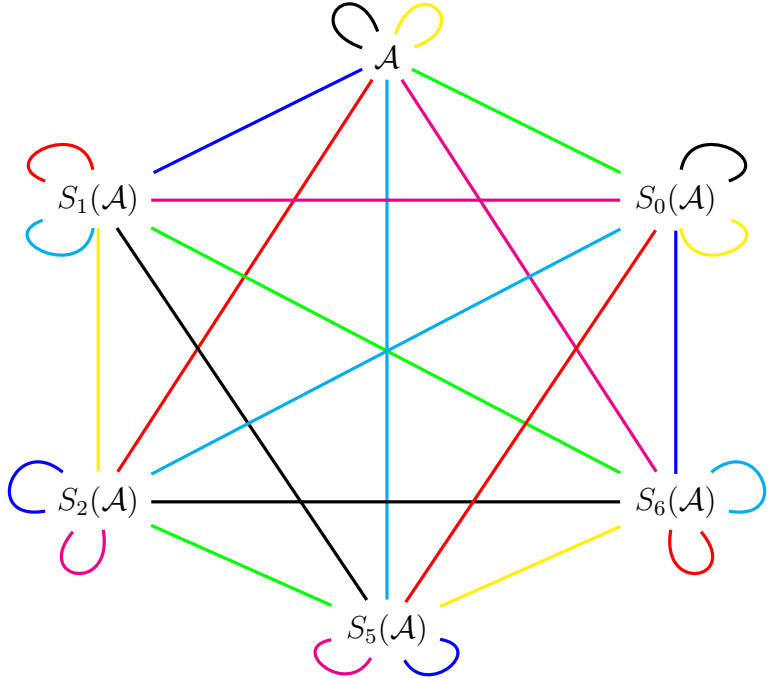 
  
 By applying the simple facts about the symmetries stated in the last part of Section \ref{SecSym}, much more information can be deduced from this proposition. In addition to $\mathcal{A}$, the mutually disjoint sets $S_0(\mathcal{A})$, $S_1(\mathcal{A})$, $S_2(\mathcal{A})$, $S_5(\mathcal{A})$ and $S_6(\mathcal{A})$ are all asymmetric invariant sets. Given any symmetry $S \in S_G$, it is easy to check whether or not each of the asymmetric sets just mentioned breaks $S$. This can be visualized with a graph, where the vertices are the sets $\mathcal{A}$, $S_0(\mathcal{A})$, $S_1(\mathcal{A})$, $S_2(\mathcal{A})$, $S_5(\mathcal{A})$ and $S_6(\mathcal{A})$, while the edges indicate the symmetries mapping the set on one end to the set on the other end. The graph is depicted on Figure \ref{FigAsym}. 
 
 Now consider a set $\mathcal{B} \in \{ \mathcal{A}, S_0(\mathcal{A}), S_1(\mathcal{A}), S_2(\mathcal{A}) , S_5(\mathcal{A}),S_6(\mathcal{A}) \}$ and a symmetry $S \in S_G$. Write $S=S_{i_j}\dots S_{i_1}$, where $S_{i_k} \in \{S_0\dots S_6\},$ $k=1\dots j$. Starting from the node $\mathcal{B}$ in our graph following the edges with labels $S_{i_k}$, $k=1\dots j$ we obtain a path. The symmetry $S$ is not broken by $\mathcal{B}$ if and only if this path ends in $\mathcal{B}$.
 
 It is also an interesting fact, that for any of the asymmetric sets the symmetry subgroup that leaves this set invariant is an abelian subgroup of $S_G$ generated by two symmetries originating from commuting permutations of the original coordinates $x_1,x_2,x_3$ and $x_4$. 
 
 \end{rem}
 
 The proof of the proposition can be found in Appendix B. Notice that Proposition 1 implies Theorem 1. 
 
 We finish this section by stating a conjecture concerning the appearance of further asymmetric invariant sets at some higher value of the coupling parameter. 
 
 It is easy to see that faces of the cube are invariant with respect to the dynamics. Explicit calculations of the \emph{2-dimensional dynamics} of these subsets show that when
 \[
 \varepsilon=\frac{5-\sqrt{17}}{2},
  \]
 new asymmetric invariant sets (polygons) appear. Regarding the 3-dimensional system, this might mean (by continuity of the dynamics) that positive-Lebesgue measure asymmetric invariant sets also appear at this value of $\varepsilon$. By our simulations it seems likely that if these polyhedra exist, the previously mentioned polygons are actually their faces.  
 
 \begin{figure}[h!!]
    \centering
    \begin{tikzpicture}[>=stealth',shorten >=1pt,auto,node distance=4cm,
                    very thick,main node/.style={circle,draw,font=\L_{\varepsilon}arge\bfseries},scale=0.9]
    
      \node[draw=none,fill=none] (00) { $\mathcal{A}_2$};
      \node[draw=none,fill=none] (0) at  (0,-8.5) {$S_0(\mathcal{A}_2)$};
      \node[draw=none,fill=none] (3) at (4,-1.5) {$S_3(\mathcal{A}_2)$};
      \node[draw=none,fill=none] (60) at (-4,-1.5) {$S_6S_0(\mathcal{A}_2)$};
      \node[draw=none,fill=none] (50) at  (5,-4) {$S_5S_0(\mathcal{A}_2)$};
      \node[draw=none,fill=none] (30) at (-4,-6.5) {$S_3S_0(\mathcal{A}_2)$};
      \node[draw=none,fill=none] (6) at   (4,-6.5) {$S_6(\mathcal{A}_2)$};
      \node[draw=none,fill=none] (5) at  (-5,-4) {$S_5(\mathcal{A}_2)$};

      \path
        (00) edge [blue, out=40,in=0,looseness=8] node {} (00)
        (00) edge [red, out=110,in=70,looseness=10] node {} (00)
        (00) edge [yellow, out=180,in=140,looseness=8] node {} (00)
        (00) edge [green] node {} (0)
        (00) edge [magenta] node {} (6)
        (00) edge [cyan] node {} (5)
        (00) edge [black] node {} (3)
        
        (60) edge [black, out=230,in=190,looseness=3] node {} (60)
             (60) edge [blue, out=120,in=160,looseness=3] node {} (60)
             (60) edge [cyan, out=50,in=90,looseness=6] node {} (60)
             (60) edge [green] node {} (6)
             (60) edge [red] node {} (30)
             (60) edge [magenta] node {} (0)
             (60) edge [yellow] node {} (50)
             
             (0) edge [yellow, out=180,in=220,looseness=4] node {} (0)
                       (0) edge [red, out=250,in=290,looseness=7] node {} (0)
                       (0) edge [blue, out=320,in=0,looseness=4] node {} (0)
                       (0) edge [cyan] node {} (50)
                       (0) edge [black] node {} (30)
                       
         (3) edge [yellow, out=130,in=90,looseness=6] node {} (3)
                             (3) edge [cyan, out=60,in=20,looseness=4] node {} (3)
                             (3) edge [magenta, out=350,in=310,looseness=4] node {} (3)
                            (3) edge [green] node {} (30)
                             (3) edge [blue] node {} (5)
                             (3) edge [red] node {} (6)
                             
        (30) edge [yellow, out=180,in=140,looseness=3] node {} (30)
                                  (30) edge [cyan, out=250,in=210,looseness=5] node {} (30)
                                  (30) edge [magenta, out=320,in=280,looseness=6] node {} (30)
                                 (30) edge [blue] node {} (50)
     
     (50) edge [magenta, out=90,in=50,looseness=6] node {} (50)
                                    (50) edge [red, out=340,in=20,looseness=3] node {} (50)
                                    (50) edge [black, out=270,in=310,looseness=6] node {} (50)
                                   (50) edge [green] node {} (5)
     
     (5) edge [magenta, out=90,in=130,looseness=6] node {} (5)
                                      (5) edge [red, out=160,in=200,looseness=3] node {} (5)
                                      (5) edge [black, out=270,in=230,looseness=6] node {} (5)
                                     (5) edge [yellow] node {} (6)
     
     (6) edge [black, out=220,in=260,looseness=6] node {} (6)
                         (6) edge [blue, out=290,in=330,looseness=4] node {} (6)
                         (6) edge [cyan, out=0,in=40,looseness=4] node {} (6);

    \end{tikzpicture}
    \caption{Eight conjectured asymmetric invariant sets and the symmetries connecting them. Edge colors indicate the following symmetries: green -- $S_0$, blue -- $S_1$, red -- $S_2$, black -- $S_3$, yellow -- $S_4$, cyan -- $S_5$, magenta -- $S_6$.} \label{FigAsym2}
    \end{figure}
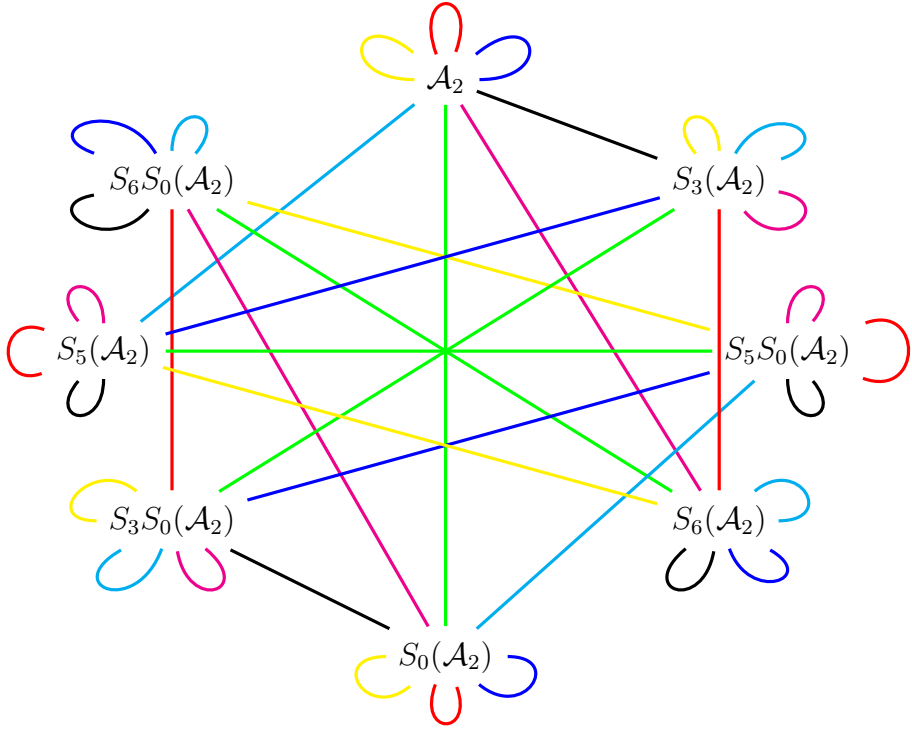 
 
\begin{con} \label{conestar2}
 There exists a polyhedron $P_3$, such that
 \[
 \mathcal{A}_2=P_3 \cup S_1(P_3) \cup S_2(P_3) \cup S_4(P_3) \cup S_2S_1(P_3) \cup S_1S_2(P_3)
 \]
 is symmetric with respect to $S_1,S_2$ and $S_4$, asymmetric with respect to $S_0,S_3,S_5$ and $S_6$, and invariant if
 \[
 0.438 \approx  \frac{5-\sqrt{17}}{2}=\varepsilon^{**} \leq \varepsilon. 
 \]
 \end{con}

 \begin{rem}
 Notice that the symmetry subgroup leaving $\mathcal{A}_2$ invariant is generated by  $S_1$ and $S_2$, which correspond to the permutations of the original coordinates $(x_1,x_2,x_3,x_4)$ such that $x_4$ is fixed. In conclusion, the symmetry subgroups leaving the symmetric images of $\mathcal{A}_2$ invariant also correspond to permutations that leave one of the coordinates $x_1,x_2,x_3$ or $x_4$ fixed.
 \end{rem}
 
 \begin{rem}
 The conjecture implies the existence of 7 further asymmetric invariant sets related to $\mathcal{A}_2$ by the structure of the symmetry group, as can be seen on the graph on Figure \ref{FigAsym2}.  
 \end{rem}
 
 We note that if a description of $P_3$ would be obtained, the proof of Conjecture 1 could be completed easily.
 
 \subsection{A non-trivial symmetric invariant set} \label{SecSymG}
 In this section, we are going to define a symmetric invariant set connected to Conjecture 2, which is stated below.
 
 We define this set as the union of all symmetric images of certain polyhedra with faces parallel to singularities. To ensure invariance, we choose the exact locations of these faces with the help of the intervals defined by Equation \eqref{mixing}. These are the intervals that support the invariant measure of the map $L_{\varepsilon}$ when $1-\frac{\sqrt{2}}{2} \leq \varepsilon < 1-\frac{\sqrt[4]{2}}{2}$. To ensure symmetry, we construct the set as the union of a polyhedron $P_0$ and all of its symmetric images under the elements of $S_G$.
 
 Let this polyhedron $P_0$ be defined as follows: 
    
 \begin{center}
     \begin{tabular}{| l | l |}
       \hline 
       & $\mathbf{P_0}$ \\ \hline
       $\mathbf{p}$ & $L_{\varepsilon}(\varepsilon/2) < p < L_{\varepsilon}(1-\varepsilon/2)$ \\ \hline
       $\mathbf{q}$ & $\varepsilon/2 < q < L_{\varepsilon}^2(1-\varepsilon/2)$ \\ \hline
       $\mathbf{r}$ & $L_{\varepsilon}(\varepsilon/2) < r < L_{\varepsilon}(1-\varepsilon/2)$ \\ \hline
       $\mathbf{p+q+r}$  & $1+\varepsilon/2 < p+q+r < 1+L_{\varepsilon}^2(1-\varepsilon/2)$  \\ \hline
         \end{tabular}
   \end{center}
 
 \begin{prop} \label{propsim}
 The set 
 \begin{align*}
 \mathcal{S}=&P_0\cup S_0(P_0) \cup S_1(P_0) \cup S_2(P_0) \cup S_3(P_0) \cup S_4(P_0) \cup S_5(P_0)  \\
 & S_0S_1(P_0) \cup S_2S_1(P_0) \cup S_3S_1(P_0) \cup S_4S_1(P_0) \cup  S_5S_1(P_0)
 \end{align*}
 is symmetric and it is invariant with respect to $G_{\varepsilon,3}$ if and only if 
 \[
 1-\frac{\sqrt{2}}{2} \leq \varepsilon.
 \]
 \end{prop}
 
\begin{figure}
   \centering
   \includegraphics[scale=0.65]{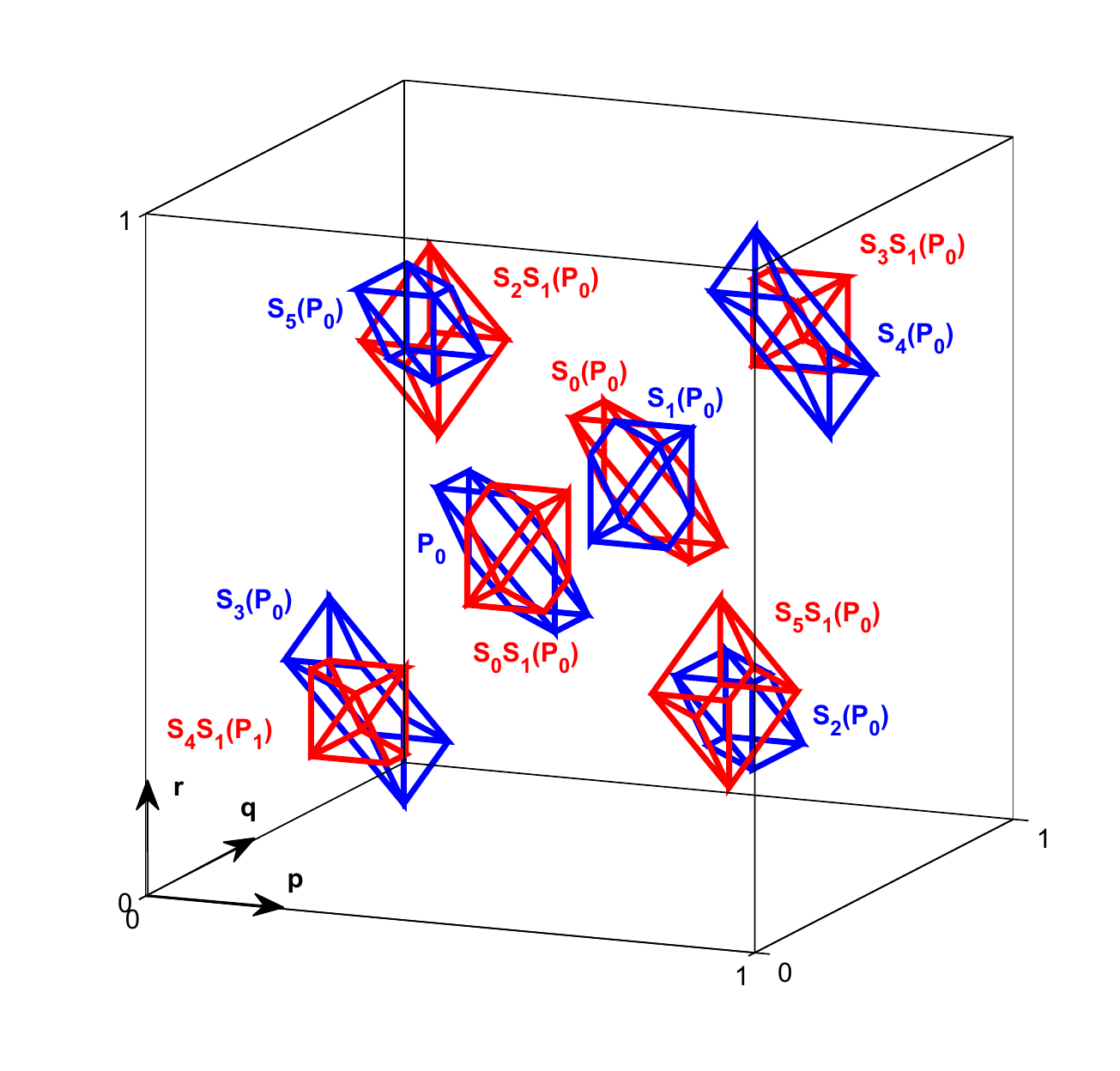}
   \caption{The symmetric set $\mathcal{S}$ for $1-\frac{\sqrt{2}}{2} \leq \varepsilon$.} \label{FigS}
   \end{figure}

 \begin{rem}
 Provided that $1-\frac{\sqrt{2}}{2} \leq \varepsilon$, the set $\mathcal{S}$ is the union of 12 convex polyhedra, see Figure \ref{FigS}. The union has 6 disjoint connected components, each of which consists of 2 intersecting convex polyhedra:
 \begin{align*}
 \text{component 1}&: \hspace{0.2cm} P_0 \cup S_0S_1(P_0), \\
 \text{component 2}&: \hspace{0.2cm} S_0(P_0) \cup S_1(P_0), \\
 \text{component 3}&: \hspace{0.2cm} S_2(P_0) \cup S_5S_1(P_0), \\ 
 \text{component 4}&: \hspace{0.2cm} S_3(P_0) \cup S_4S_1(P_0), \\
 \text{component 5}&: \hspace{0.2cm} S_4(P_0) \cup S_3S_1(P_0), \\
 \text{component 6}&: \hspace{0.2cm} S_5(P_0) \cup S_2S_1(P_0).
 \end{align*}
 \end{rem}
 
 The proof of the proposition can be found in Appendix C. 
 
  An additional motivation for describing this set is to state an interesting conjecture on the possible existence of infinitely many values of the coupling parameter, where new invariant sets appear. 
 
 We have stated in Section \ref{SecL}, that when $\varepsilon=1-\frac{\sqrt{2}}{2}$, the invariant measure of $L_{\varepsilon}$ obtains a support consisting of the union of two mixing components. The symmetric set $\mathcal{S}$, which is defined with the help of these intervals, becomes invariant at the exact same value of $\varepsilon$. It would be of great interest if one could shed some light on the complicated connection between these two phenomena. The proof in Appendix C contains some ideas that might help the investigation of this question.  
 
 We have also stated in Section \ref{SecL} that at every value 
 \[
  \varepsilon_n=1-\frac{\sqrt[2^n]{2}}{2},
 \]
  the support of the invariant measure of $L_{\varepsilon}$ becomes the union of twice as many mixing components than previously. This is most likely to be the cause of some numerically observed changes in the structure of the set $\mathcal{S}$ at these values of the coupling parameter. 
  
  To enlighten this somewhat, we ask the reader to consider a three-dimensional system with every coordinate evolving according to $L_{\varepsilon}$. In a system like this, new (symmetric) invariant sets appear at every $\varepsilon_n$ due to the appearance of new mixing components of the coordinate maps. Now on each component of $\mathcal{S}$, if we apply a change of coordinates (depending on the component), $G_{\varepsilon,3}$ acts as the map previously described. But the connection of the two systems is not clear and proved too complicated to explore. However, our simulations indicate that similar phenomena take place.
 
 \begin{con}At every
 \[
  \varepsilon_n=1-\frac{\sqrt[2^n]{2}}{2} \quad n>1,
 \]
  new symmetric invariant sets appear. These new invariant sets are contained in $\mathcal{S}$ and they are symmetric. 
 \end{con}

\section{Concluding remarks}

The goal of this paper was to contribute to the understanding of a globally coupled map with four sites, a special case of the model introduced by Koiller and Young in \cite{koiller2010coupled}. In particular, our aim was to show that for large values of the coupling parameter (indicating strong interaction between the sites), asymmetric invariant sets exist. These sets imply that the invariant measure of maximal support, which is well known to be ergodic and mixing for weak coupling, is in fact not ergodic for sufficiently strong coupling. This phenomenon was not present for the model with two sites, but the coupling of three sites did produce it. Proving ergodicity breaking for four sites assures us that it is not an artifact of a somewhat low-dimensional system, but indeed a possible tendency.

In this paper, we described a value $\varepsilon^* \approx 0.397$, such that multiple asymmetric invariant sets exist when $\varepsilon \geq \varepsilon^*$. For such values of the coupling parameter, we gave an asymmetric invariant set $\mathcal{A}$ which has a simple geometric structure, namely it is the union of polyhedra. Finding this set was a highly nontrivial task, a great amount of simulations led to the conjecture of its precise parameters. We showed that this set $\mathcal{A}$ is not fully asymmetric, it is symmetric with respect to a special abelian subgroup of the map's symmetry group, which is generated by two commuting permutation symmetries. Comparison with the numerical results of Fernandez shows that it is very likely that we succeeded in giving the exact value of $\varepsilon_a$ (the smallest value of $\varepsilon$ where the existence of asymmetric invariant sets were indicated by simulations). However, it remains to be shown that no asymmetric invariant sets exist for $\varepsilon < \varepsilon^*$. It would be also interesting to obtain more insight on the role of the Lorenz map which was used in the construction of $\mathcal{A}$.

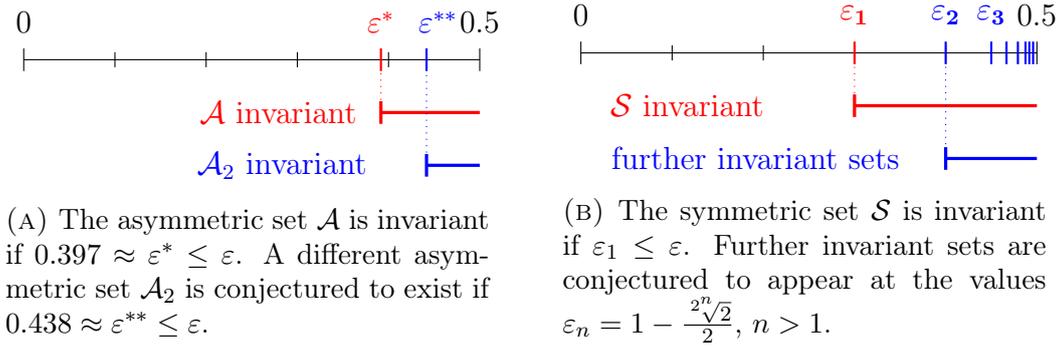
\begin{figure}[h!]
\centering
\begin{subfigure}[b]{0.4\textwidth}
\begin{tikzpicture}
\draw (0,0) -- (6,0);
\draw (0,-0.15) -- (0,0.15);
\draw (6,-0.15) -- (6,0.15);
\draw (6/5,-0.1) -- (6/5,0.1);
\draw (2*6/5,-0.1) -- (2*6/5,0.1);
\draw (3*6/5,-0.1) -- (3*6/5,0.1);
\draw (4*6/5,-0.1) -- (4*6/5,0.1);
\draw (0,0.5) node {0};
\draw (6,0.5) node {0.5};
\draw[red] (4.7,0.5) node {$\mathbf{\varepsilon^*}$};
\draw[thick, red] (4.7,-0.15) -- (4.7,0.15) ;
\draw[blue] (5.45,0.5) node {$\mathbf{\varepsilon^{**}}$};
\draw[thick, blue] (5.3,-0.15) -- (5.3,0.15) ;
\draw[very thick, red] (4.7,-0.85) -- (4.7,0.15-0.7) ;
\draw[very thick, red] (4.7,-0.7) -- (6,-0.7);
\draw[red] (3.35,-0.7) node {$\mathcal{A}$ invariant};
\draw[dotted,red] (4.7,-0.7) -- (4.7,-0.1);
\draw[very thick, blue] (5.3,-0.85-0.7) -- (5.3,0.15-1.4) ;
\draw[very thick, blue] (5.3,-1.4) -- (6,-1.4);
\draw[blue] (3.4,-1.4) node {$\mathcal{A}_2$ invariant};
\draw[dotted,blue] (5.3,-1.4) -- (5.3,-0.1);
\end{tikzpicture}
\caption{The asymmetric set $\mathcal{A}$ is invariant if $0.397 \approx \varepsilon^* \leq \varepsilon$. A different asymmetric set $\mathcal{A}_2$ is conjectured to exist if $0.438 \approx \varepsilon^{**} \leq \varepsilon$. }
\end{subfigure}
\qquad
\begin{subfigure}[b]{0.4\textwidth}
\begin{tikzpicture}
\draw (0,0) -- (6,0);
\draw (0,-0.15) -- (0,0.15);
\draw (6,-0.15) -- (6,0.15);
\draw (6/5,-0.1) -- (6/5,0.1);
\draw (2*6/5,-0.1) -- (2*6/5,0.1);
\draw (3*6/5,-0.1) -- (3*6/5,0.1);
\draw (4*6/5,-0.1) -- (4*6/5,0.1);
\draw (0,0.5) node {0};
\draw (6,0.5) node {0.5};
\draw[red] (3.6,0.5) node {$\mathbf{\varepsilon_1}$};
\draw[thick, red] (3.6,-0.15) -- (3.6,0.15) ;
\draw[very thick, red] (3.6,-0.85) -- (3.6,0.15-0.7) ;
\draw[very thick, red] (3.6,-0.7) -- (6,-0.7);
\draw[dotted, red] (3.6,-0.1) -- (3.6,-0.7);
\draw[dotted, blue] (4.8,-0.1) -- (4.8,-1.4);
\draw[red] (1.4,-0.7) node {$\mathcal{S}$ invariant};
\draw[blue] (4.8,0.5) node {$\mathbf{\varepsilon_2}$};
\draw[thick, blue] (4.8,-0.15) -- (4.8,0.15) ;
\draw[very thick, blue] (4.8,-0.85-0.7) -- (4.8,0.15-1.4) ;
\draw[very thick, blue] (4.8,-1.4) -- (6,-1.4);
\draw[blue] (2.3,-1.4) node {further invariant sets};
\draw[blue] (5.4,0.5) node {$\mathbf{\varepsilon_3}$};
\draw[thick, blue] (5.4,-0.15) -- (5.4,0.15) ;
\draw[thick, blue] (5.6,-0.15) -- (5.6,0.15) ;
\draw[thick, blue] (5.75,-0.15) -- (5.75,0.15) ;
\draw[thick, blue] (5.85,-0.15) -- (5.85,0.15) ;
\draw[thick, blue] (5.9,-0.15) -- (5.9,0.15) ;
\draw[thick, blue] (5.95,-0.15) -- (5.95,0.15) ;
\end{tikzpicture}
\caption{The symmetric set $\mathcal{S}$ is invariant if $\varepsilon_1 \leq \varepsilon$. Further invariant sets are conjectured to appear at the values $\varepsilon_n=1-\frac{\sqrt[2^n]{2}}{2}$, $n > 1$.}
\end{subfigure}
\caption{Overview of the critical parameters and the corresponding invariant sets. Results are marked with red, conjectures with blue.} \label{FigOver}
\end{figure}

Based on explicit calculations on certain lower dimensional invariant subsets of the phase space, we stated a very particular conjecture for the second critical value of the coupling parameter where new asymmetric invariant sets emerge. Our conjectured value is in good agreement with $\varepsilon_b$ (the value of $\varepsilon$ where the existence of new asymmetric invariant sets were indicated by the simulations of Fernandez). 


We also described a nontrivial symmetric invariant set $\mathcal{S}$ derived from an underlying one-dimensional Lorenz map. The existence of this set seemed clear from simulating orbits of random phase points, but finding the  right parameters of this set was nontrivial. Based on simulations, we conjecture that new invariant sets arise inside this set for countably many values of the coupling parameter. A better understanding of how the Lorenz map shapes the dynamics of this symmetric invariant set is needed to prove this particularly interesting conjecture.

To summarize, the following picture seems likely: for $\varepsilon_1=1-\sqrt{2}/2$, the attractor is exactly the set $\mathcal{S}$ (for lower values of the coupling parameter it is much larger). When $\varepsilon_1 \leq \varepsilon < \varepsilon^* \approx 0.397$, the attractor is contained in $\mathcal{S}$. The system is transitive on the attractor in all of these cases. We proved that at the value $\varepsilon^*$ the set $\mathcal{A}$ becomes invariant. By straightforward calculations we can obtain that $\mathcal{A}$ and $\mathcal{S}$ are disjoint. Hence new invariant sets emerge in locations where no trajectory returned for $\varepsilon < \varepsilon^*$. At the value of $\varepsilon^{**} \approx 0.438$, similar phenomena can be observed as in the case of $\varepsilon^*$: new asymmetric invariant sets appear in locations left by all trajectories for $\varepsilon < \varepsilon^{**}$. Furthermore, eventually the symmetric invariant set also decomposes into smaller, but still symmetric invariant sets. Some of these statements remain to be proved. A summary of the existing results and conjectures is depicted in Figure \ref{FigOver}. 
 
We would like to point out that studying this system provided us information on what phenomena seen in the system of three sites was the specialty of that system and what might be more general. By simulations of this system it is clear that \emph{ergodicity does not necessarily break down in one step} (as it happened in the $N=3$ case), but further invariant sets can emerge for higher values of the coupling parameter than the bifurcation value $\varepsilon^*$ (as stated in Conjecture 1.). We have seen that the \emph{asymmetric invariant sets do not necessarily arise as the decomposition of the symmetric invariant set}, as seen in the case of three sites, but they might appear in a part of the phase space which was transient for $1-\frac{\sqrt{2}}{2} \leq \varepsilon < \varepsilon^*$. A further interesting observation (by simulations) is that the symmetric invariant set decomposes into several, although symmetric invariant sets. 

Some of the interesting phenomena observed in simulations remain as conjectures, so although a proof of ergodicity breaking was obtained, there is possibility for future work even in this particular system.

\section*{Acknowledgements}
The author is grateful for motivating discussions with Bastien Fernandez and P\'eter B\'alint. The author would also like to express gratitude towards an anonymous referee whose many helpful comments improved the paper substantially. This work was partially supported by
Hungarian National Foundation for Scientific Research (NKFIH OTKA) grants K104745 \& K123782, and Stiftung Aktion \"Osterreich Ungarn (A\"OU) grants 87\"ou6 \& 92\"ou6.

\newpage

\bibliographystyle{plain}
\bibliography{references_n4}
 \nocite{*}
 
\newpage

\section*{Appendix A: Continuity domains of the map $G_{\varepsilon,3}$}

In this section we describe the continuity domains of $G_{\varepsilon,3}$ and the action of the map on each of them. This section provides essential reference for Appendices B and C. 

\begin{figure}[htbp]
\centering
\includegraphics[scale=1]{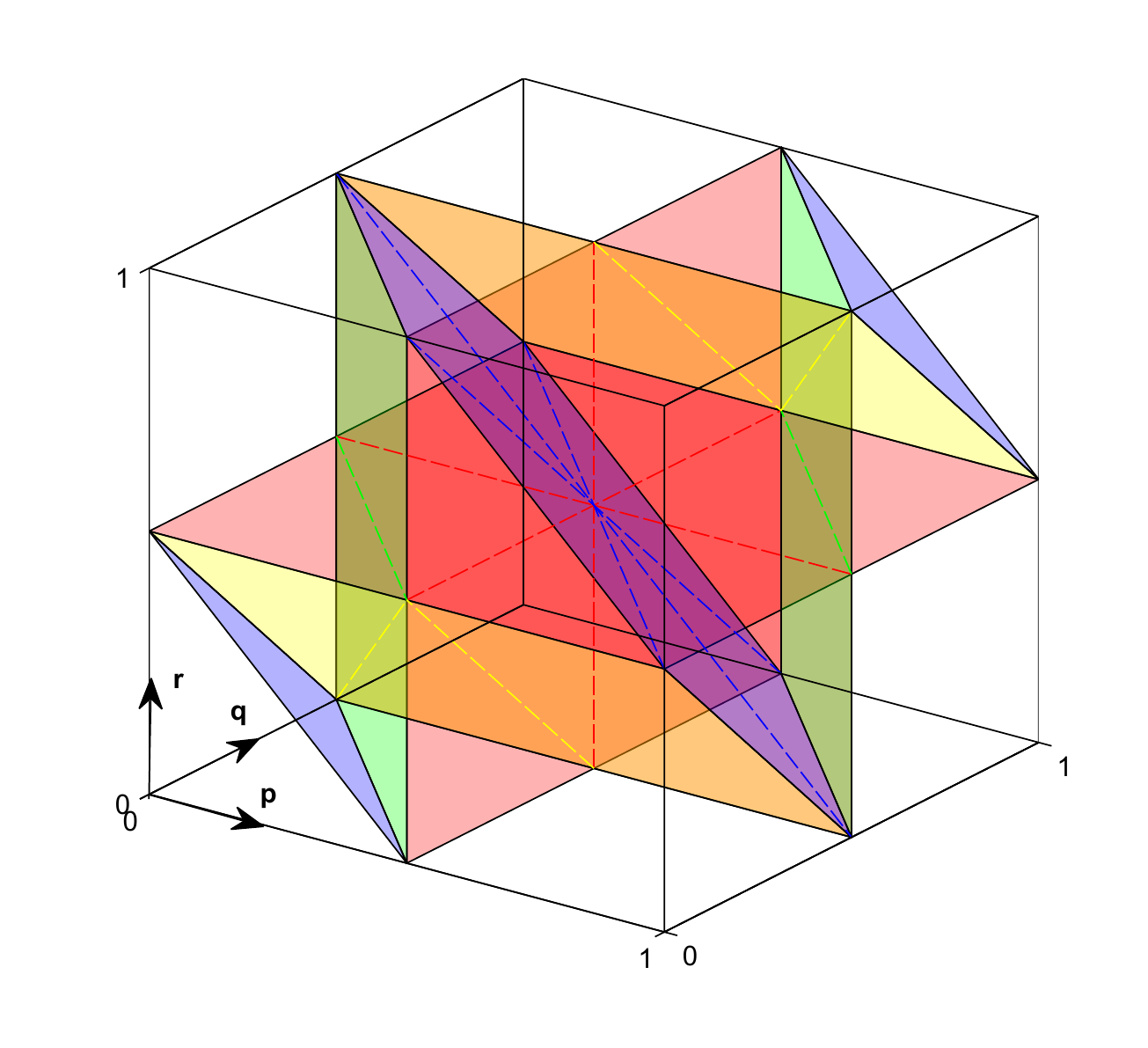}
\caption{Continuity domains of the map $G_{\varepsilon,3}$. $\mathbb{T}^3$ is represented as the unit cube of $\mathbb{R}^3$. The singularities are the intersections of the following planes with the unit cube: $p=1/2$, $q=1/2$, $r=1/2$ (red), $p+q=1/2$, $p+q=3/2$ (green), $q+r=1/2$, $q+r=3/2$ (yellow), $p+q+r=1/2$, $p+q+r=3/2$, $p+q+r=5/2$ (blue).} \label{FigSing}
\end{figure}

The map $G_{\varepsilon,3}$ is a piecewise affine map on $\mathbb{T}^3$ (represented as the unit cube of $\mathbb{R}^3$, with opposite faces identified). The singularities, arising from the function $g$, are intersections of certain planes with the unit cube. This is pictured on Figure \ref{FigSing}.

\begin{figure}[htbp]
      \centering
      \begin{subfigure}[b]{0.3\textwidth}
      \includegraphics[scale=0.55]{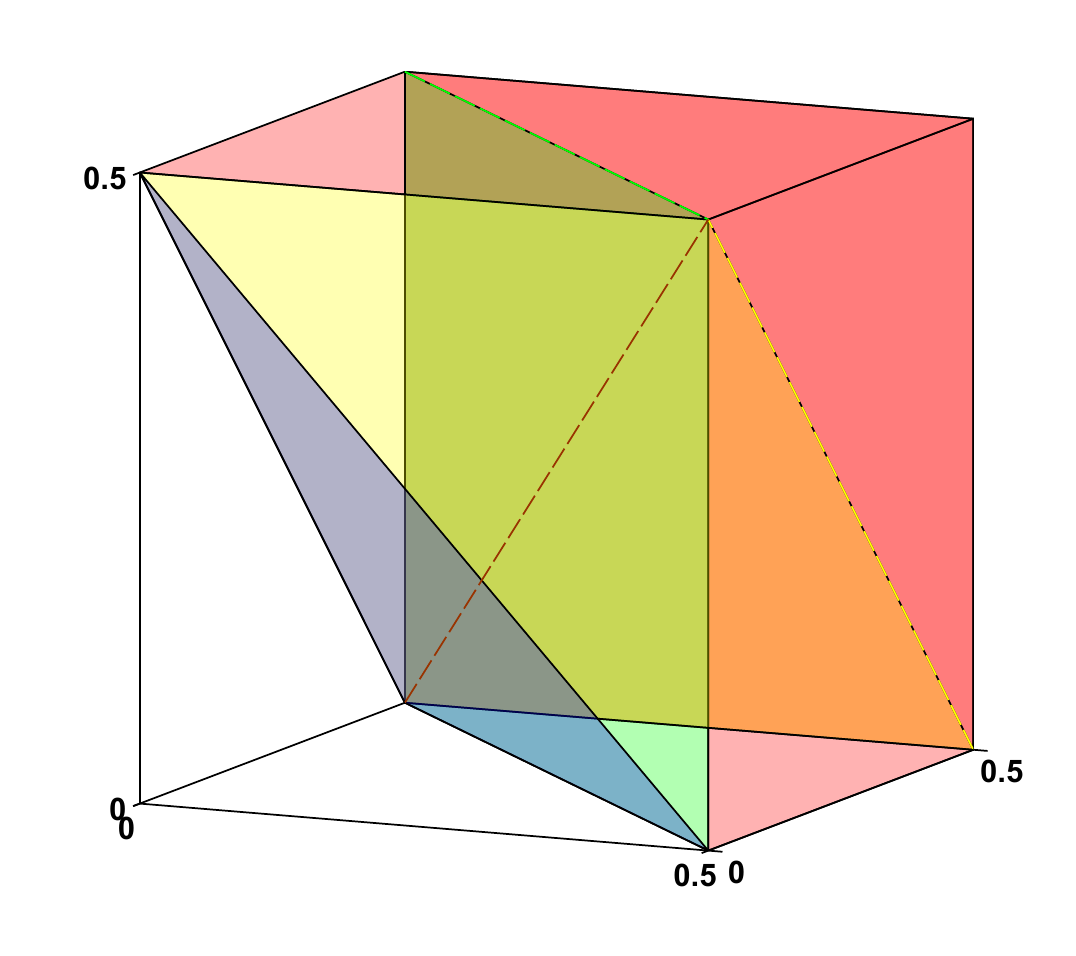}
      \caption{Cube $1$, where $0 < p,q,r < 1/2$.}
      \end{subfigure}
      \hspace{2cm}
      \begin{subfigure}[b]{0.3\textwidth}
       \includegraphics[scale=0.55]{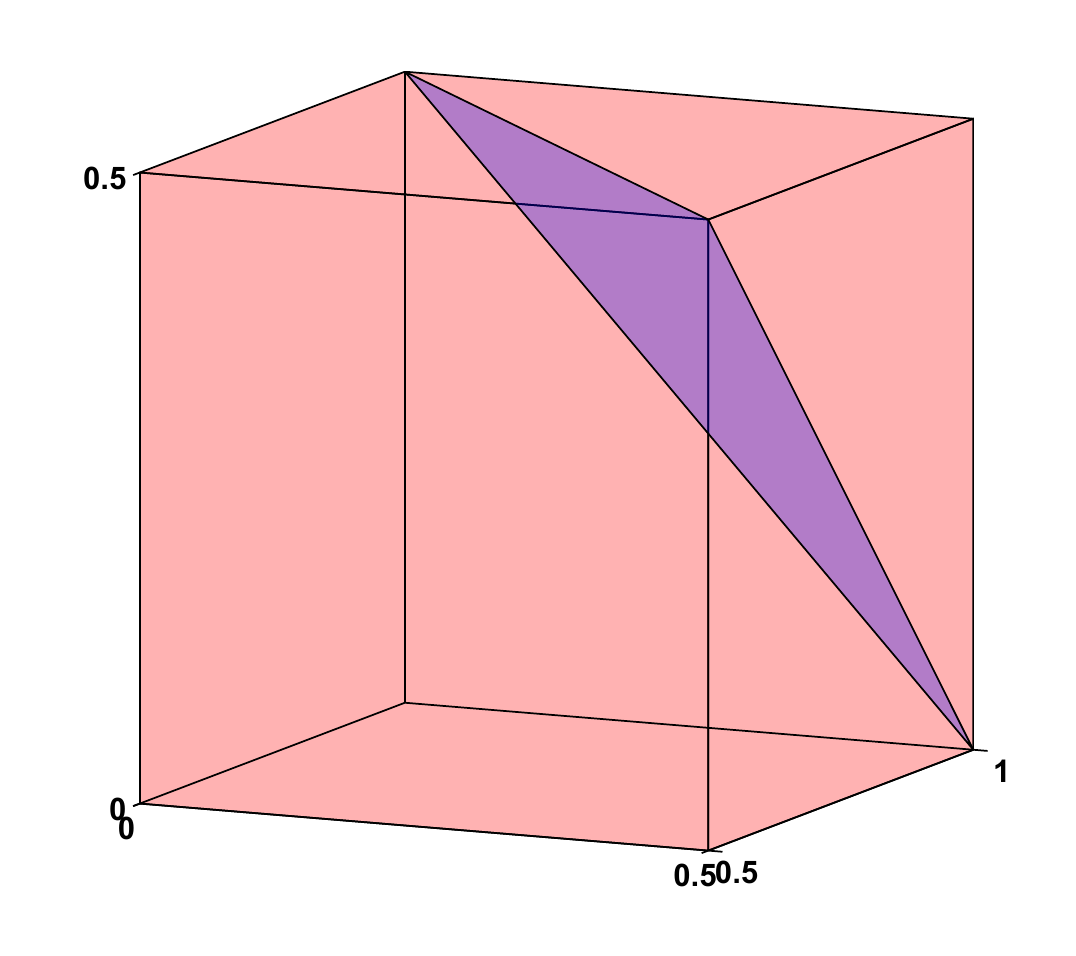}
       \caption{Cube $2$, where $0 < p,r < 1/2$ and $1/2 < q < 1$.}
       \end{subfigure}
      \caption{Cubes $1$ and $2$.} \label{FigSing12}
      \end{figure}

We are going to fix a notation for the continuity domains of the map $G_{\varepsilon,3}$ for further reference. To give these polyhedra systematically, we decompose the cube into eight smaller cubes according to the singularities $p=1/2$, $q=1/2$ and $r=1/2$ (pictured in light red on Figure \ref{FigSing}). We number each cube, and we further decompose them according to the remaining singularities, marking each final domain with a letter. Geometrically, this is pictured on Figures \ref{FigSing12}-\ref{FigSing78}. To give a precise description, each polyhedral domain can be characterized by a set of inequalities. The minimal such descriptions (in $\mathbb{R}^3$) are found in the Tables \ref{Tab12}-\ref{Tab78}.

\begin{table}[htbp]
  \begin{tabular}{| p{0.5cm} | p{2.3cm} | p{2.3cm} | p{2.5cm} | p{2cm} | p{2cm} | p{2.5cm} |}
    \hline 
     & $\mathbf{p}$ & $\mathbf{q}$ & $\mathbf{r}$ & $\mathbf{p+q}$ & $\mathbf{q+r}$ & $\mathbf{p+q+r}$ \\ \hline \hline
     $\mathbf{1a}$ & $p > 0$ & $q > 0$ & $r > 0$ &  &  & $p+q+r < 1/2$ \\ \hline
     $\mathbf{1b}$ & $p < 1/2$ & & $r > 0$ & $p+q > 1/2$ & $q+r < 1/2$ & \\ \hline
     $\mathbf{1c}$ & $p > 0$ & & $r < 1/2$ & $p+q < 1/2$ & $q+r > 1/2$ & \\ \hline
     $\mathbf{1d}$ &  & $q > 0 $ &  & $p+q < 1/2$ & $q+r < 1/2$ & $p+q+r > 1/2$\\ \hline
     $\mathbf{1e}$ & $p < 1/2$  & $q < 1/2 $ & $r < 1/2$  & $p+q > 1/2$ & $q+r > 1/2$ &\\ \hline \hline
     $\mathbf{2a}$ & $0 < p < 1/2$  & $1/2 < q < 1 $ & $0 < r < 1/2$  &  & & $p+q+r < 3/2$\\ \hline
     $\mathbf{2b}$ & $p < 1/2$  & $ q < 1 $ & $ r < 1/2$  &  & & $p+q+r > 3/2$\\ \hline 
     \end{tabular}
     \caption{Domains of continuity contained in cubes 1 and 2.} \label{Tab12}
     \end{table}
     
     \begin{table}[htbp]
       \begin{tabular}{| p{0.5cm} | p{2.3cm} | p{2.3cm} | p{2.5cm} | p{2cm} | p{2cm} | p{2.5cm} |}
         \hline 
          & $\mathbf{p}$ & $\mathbf{q}$ & $\mathbf{r}$ & $\mathbf{p+q}$ & $\mathbf{q+r}$ & $\mathbf{p+q+r}$ \\ \hline \hline
     $\mathbf{3a}$ & $p >1/2$  & $ q > 1/2 $ & $ r < 1/2$  & $p+q < 3/2$ & & $p+q+r > 3/2$ \\ \hline
     $\mathbf{3b}$ & $p > 1/2$  & $ q > 1/2 $ & $ r > 0$  &  & & $p+q+r < 3/2$\\ \hline
     $\mathbf{3c}$ & $p < 1$  & $ q < 1 $ & $ 0 < r < 1/2$  & $p+q > 3/2$ & & \\ \hline \hline
     $\mathbf{4a}$ & $1/2 < p < 1$  & $ q > 0$ & $ r > 0$  & & $q+r < 1/2$ & \\ \hline
     $\mathbf{4b}$ & $p > 1/2$  & $ q < 1/2$ & $ r < 1/2$  & & $q+r > 1/2$ & $p+q+r < 3/2$ \\ \hline
     $\mathbf{4c}$ & $p < 1$  & $ q < 1/2$ & $ r < 1/2$  & & & $p+q+r > 3/2$ \\ \hline
     \end{tabular}
     \caption{Domains of continuity contained in cubes 3 and 4.} \label{Tab34}
          \end{table}
      
      \begin{figure}[htbp]
           \centering
           \begin{subfigure}[b]{0.3\textwidth}
           \includegraphics[scale=0.55]{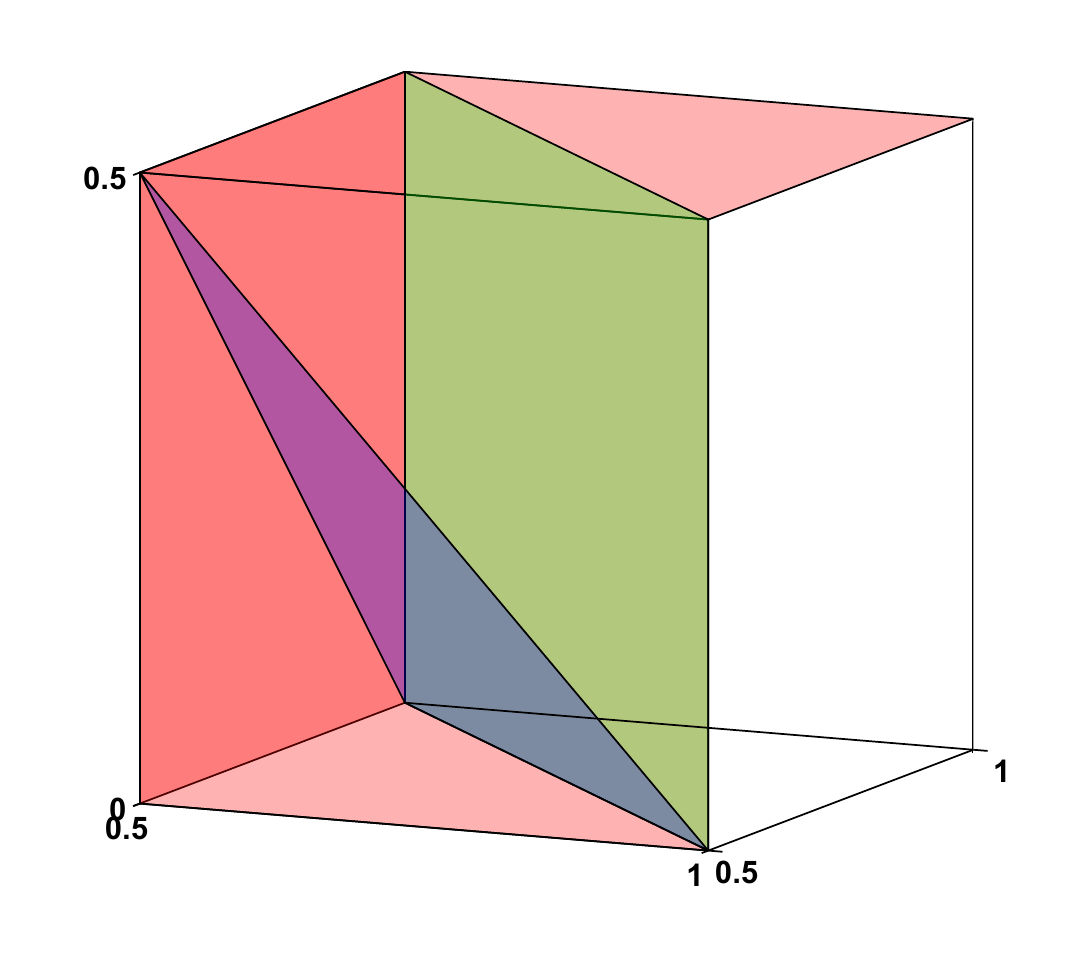}
           \caption{Cube $3$, where $0 < r < 1/2$ and $1/2 <p,q<1$.}
           \end{subfigure}
           \hspace{2cm}
           \begin{subfigure}[b]{0.3\textwidth}
            \includegraphics[scale=0.55]{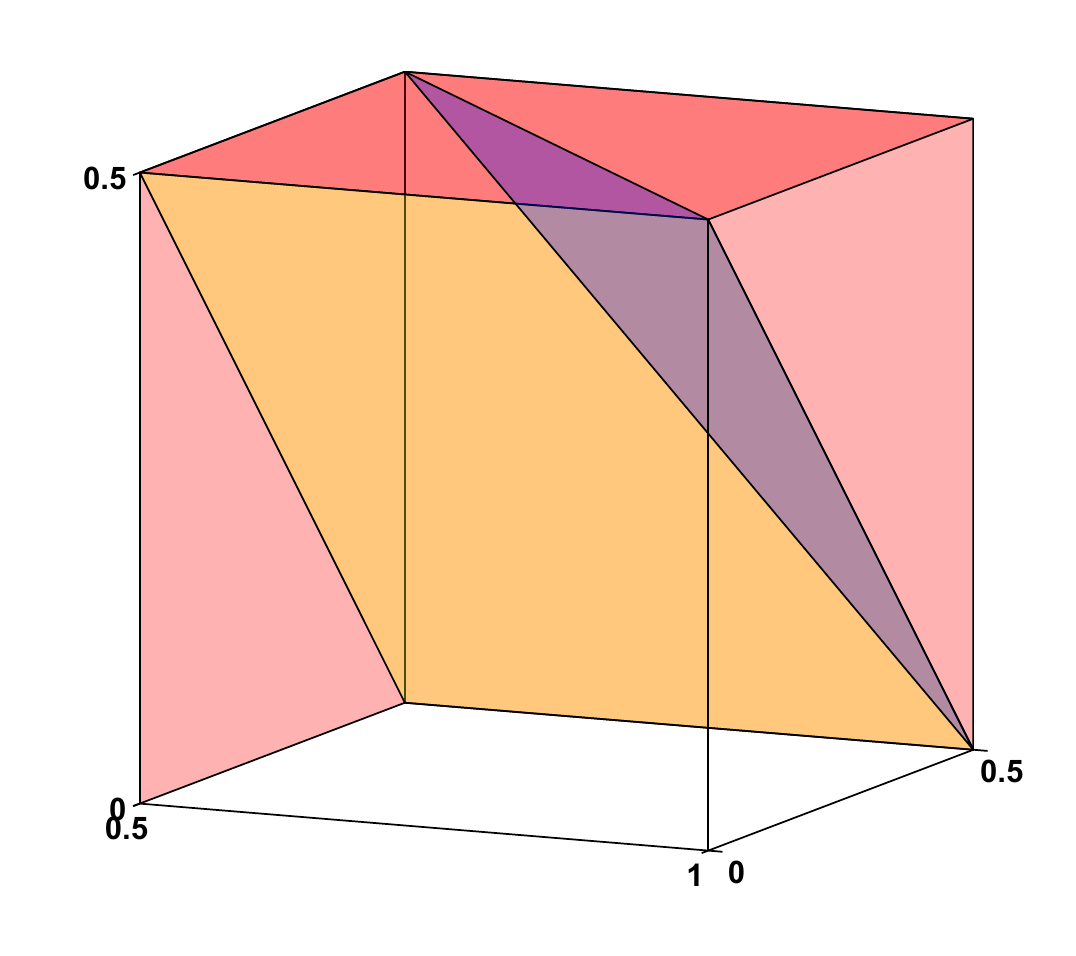}
            \caption{Cube $4$, where $0 < q,r < 1/2$ and $1/2 < p < 1$.}
            \end{subfigure}
           \caption{Cubes $3$ and $4$.}
           \end{figure}     
          
     \begin{table}[htbp]
            \begin{tabular}{| p{0.5cm} | p{2.3cm} | p{2.3cm} | p{2.5cm} | p{2cm} | p{2cm} | p{2.5cm} |}
              \hline 
               & $\mathbf{p}$ & $\mathbf{q}$ & $\mathbf{r}$ & $\mathbf{p+q}$ & $\mathbf{q+r}$ & $\mathbf{p+q+r}$ \\ \hline \hline     
     $\mathbf{5a}$ & $p > 0$  & $ q > 0$ & $1/2 < r < 1$  & $p+q < 1/2$ & & \\ \hline
     $\mathbf{5b}$ & $p < 1/2$  & $ q < 1/2$ & $r > 1/2$  & $p+q > 1/2$ & & $p+q+r < 3/2$ \\ \hline
     $\mathbf{5c}$ & $p < 1/2$  & $ q < 1/2$ & $r < 1$  & & & $p+q+r > 3/2$ \\ \hline \hline
     $\mathbf{6a}$ & $p < 1/2$  & $ q > 1/2$ & $r > 1/2$  & & $q+r < 3/2$ & $p+q+r > 3/2$  \\ \hline 
     $\mathbf{6b}$ & $p > 0$  & $ q > 1/2$ & $r > 1/2$  & & & $p+q+r < 3/2$ \\ \hline
     $\mathbf{6c}$ & $0 < p < 1/2$  & $ q < 1$ & $r < 1$  & & $q+r > 3/2$ &  \\ \hline
     \end{tabular}
     \caption{Domains of continuity contained in cubes 5 and 6.} \label{Tab56}
               \end{table}
               
     \begin{figure}[htbp]
                \centering
                \begin{subfigure}[b]{0.3\textwidth}
                \includegraphics[scale=0.55]{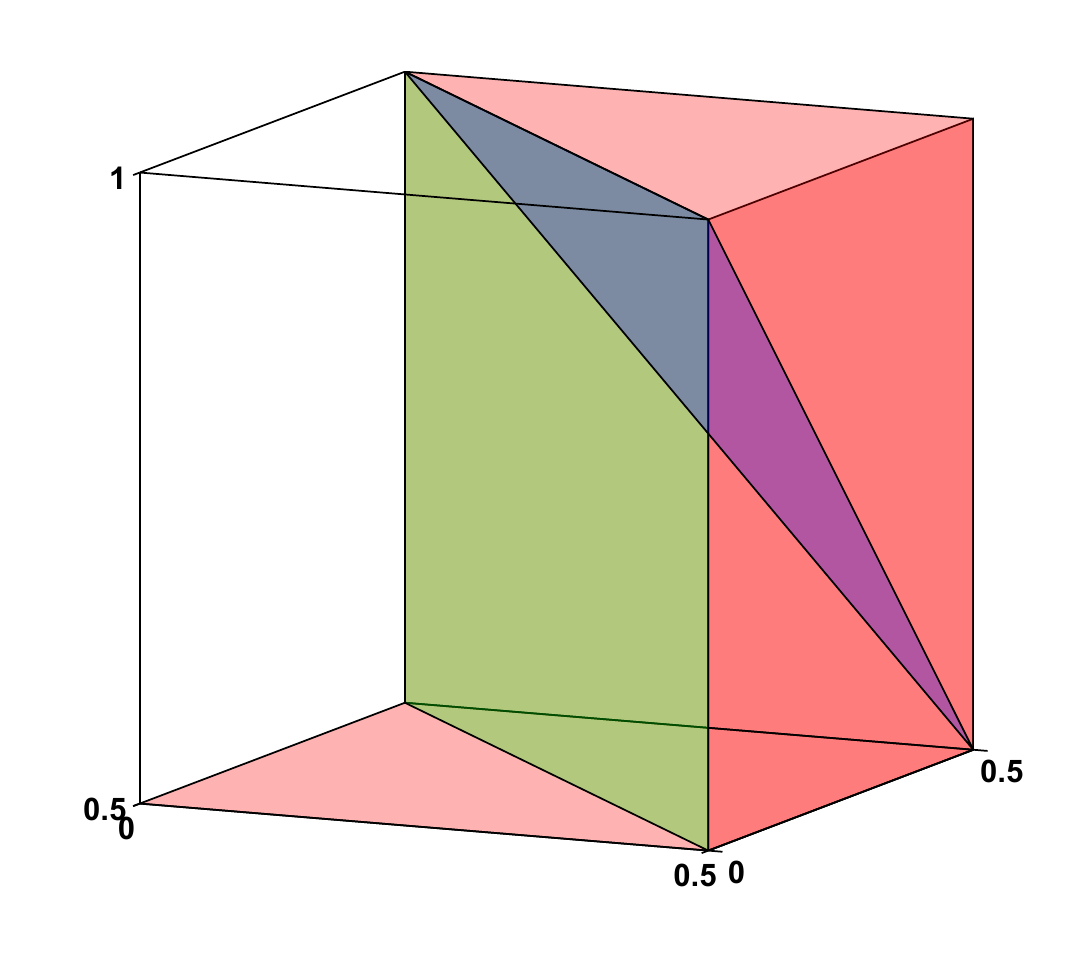}
                \caption{Cube $5$, where $0 < p,q < 1/2$ and $1/2 < r < 1$.}
                \end{subfigure}
                \hspace{2cm}
                \begin{subfigure}[b]{0.3\textwidth}
                 \includegraphics[scale=0.56]{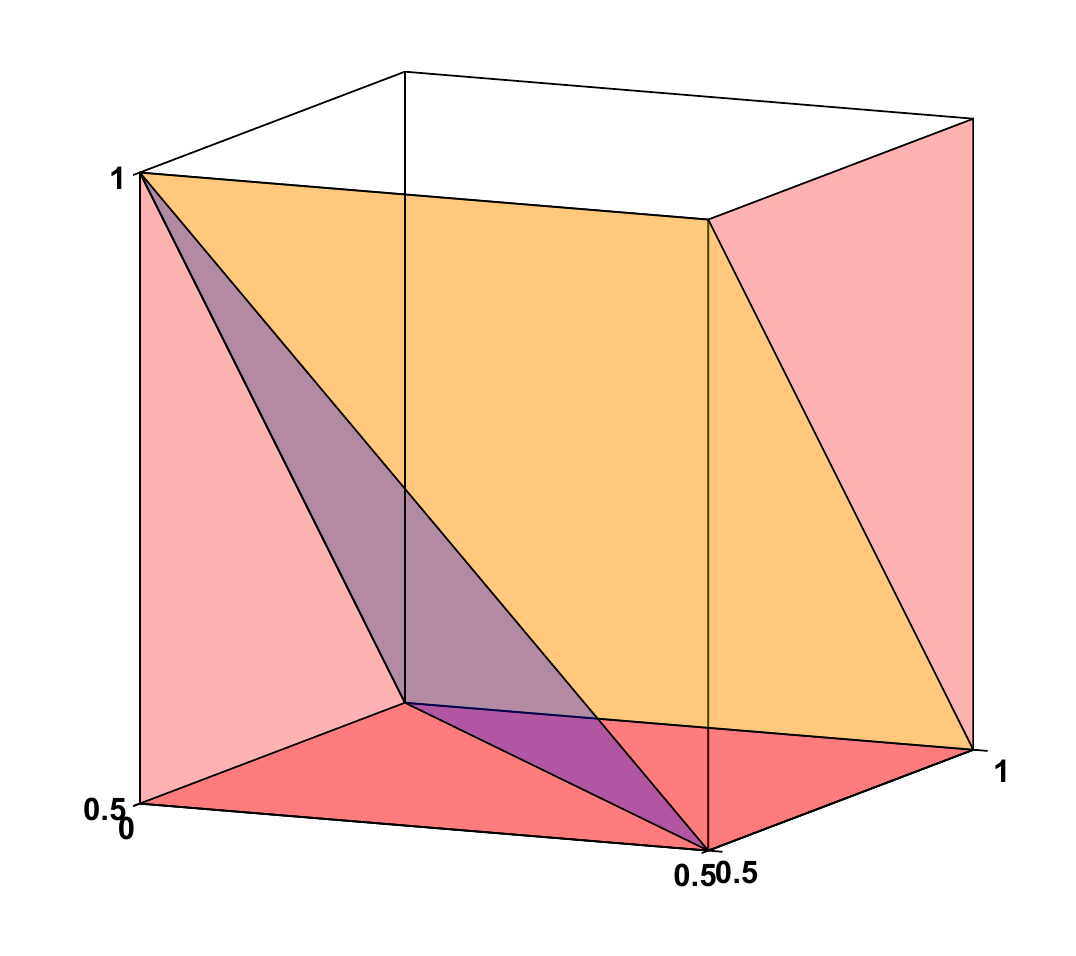}
                 \caption{Cube $6$, where $0 < p < 1/2$ and $1/2 < q,r < 1$.}
                 \end{subfigure}
                \caption{Cubes $5$ and $6$.}
                \end{figure}     
             
     \begin{table}[htbp]
                 \begin{tabular}{| p{0.5cm} | p{2.3cm} | p{2.3cm} | p{2.5cm} | p{2cm} | p{2cm} | p{2.5cm} |}
                   \hline 
                    & $\mathbf{p}$ & $\mathbf{q}$ & $\mathbf{r}$ & $\mathbf{p+q}$ & $\mathbf{q+r}$ & $\mathbf{p+q+r}$ \\ \hline \hline              
     $\mathbf{7a}$ & $p < 1$  & $ q < 1$ & $r < 1$  & & & $p+q+r > 5/2$  \\ \hline
     $\mathbf{7b}$ & $p > 1/2$  & & $r < 1$  & $p+q < 3/2$ & $q+r > 3/2$ &  \\
     \hline
     $\mathbf{7c}$ & $p < 1$  & & $r > 1/2$  & $p+q > 3/2$ & $q+r < 3/2$ &  \\ \hline
     $\mathbf{7d}$ &  & $q < 1$ &  & $p+q > 3/2$ & $q+r > 3/2$ & $p+q+r < 5/2$  \\ \hline
     $\mathbf{7e}$ & $p > 1/2$  & $q >1/2$ & $r > 1/2$  & $p+q < 3/2$ & $q+r < 3/2$ &  \\ \hline \hline
     $\mathbf{8a}$ & $1/2 < p < 1$  & $0 < q < 1$ & $1/2 < r < 1$  &  & & $p+q+r > 3/2$  \\ \hline 
     $\mathbf{8b}$ & $p > 1/2$  & $q >0$ & $r > 1/2$  &  & & $p+q+r < 3/2$  \\ \hline 
  \end{tabular}
  \caption{Domains of continuity contained in cube 7 and 8.} \label{Tab78}
\end{table}

\clearpage

\begin{figure}[htpb]
                \centering
                \begin{subfigure}[b]{0.3\textwidth}
                \includegraphics[scale=0.55]{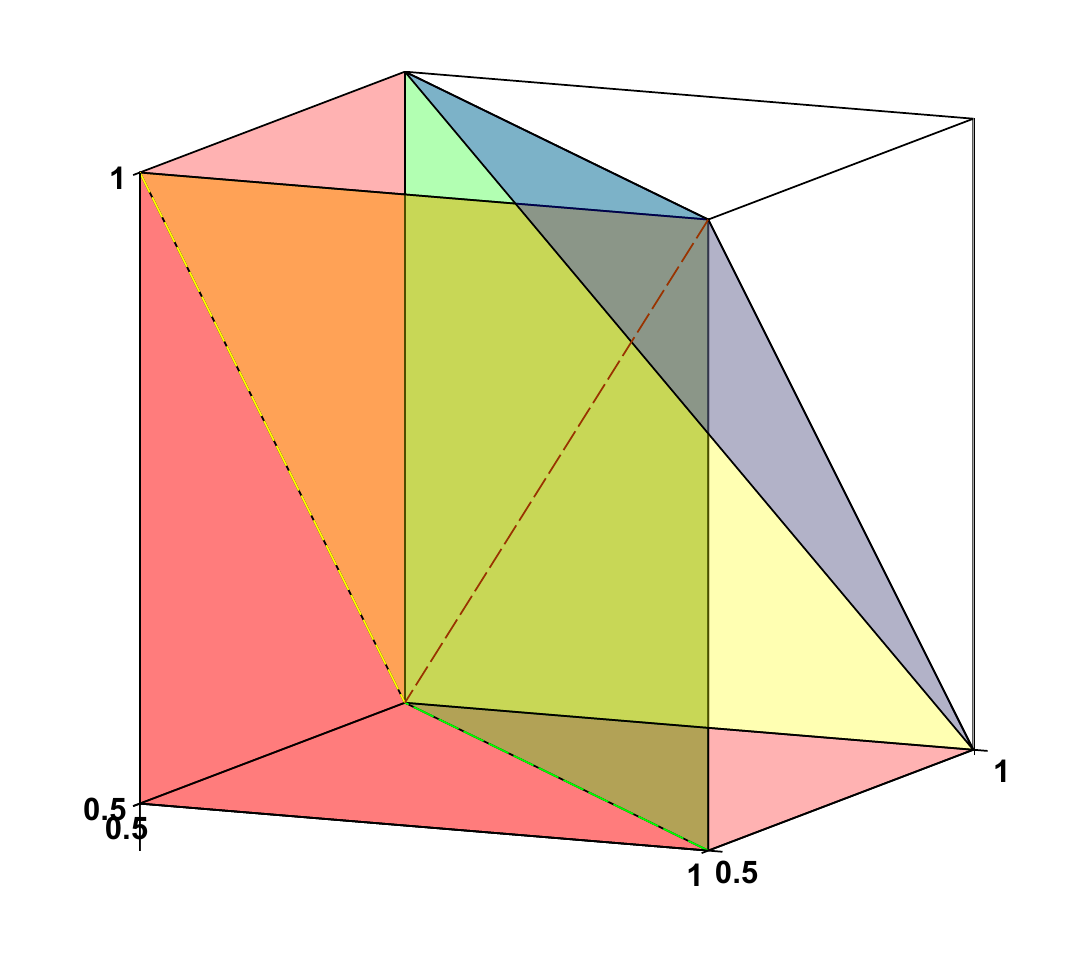}
                \caption{Cube $7$, where $1/2 < p,q,r < 1$.}
                \end{subfigure}
                \hspace{2cm}
                \begin{subfigure}[b]{0.3\textwidth}
                 \includegraphics[scale=0.55]{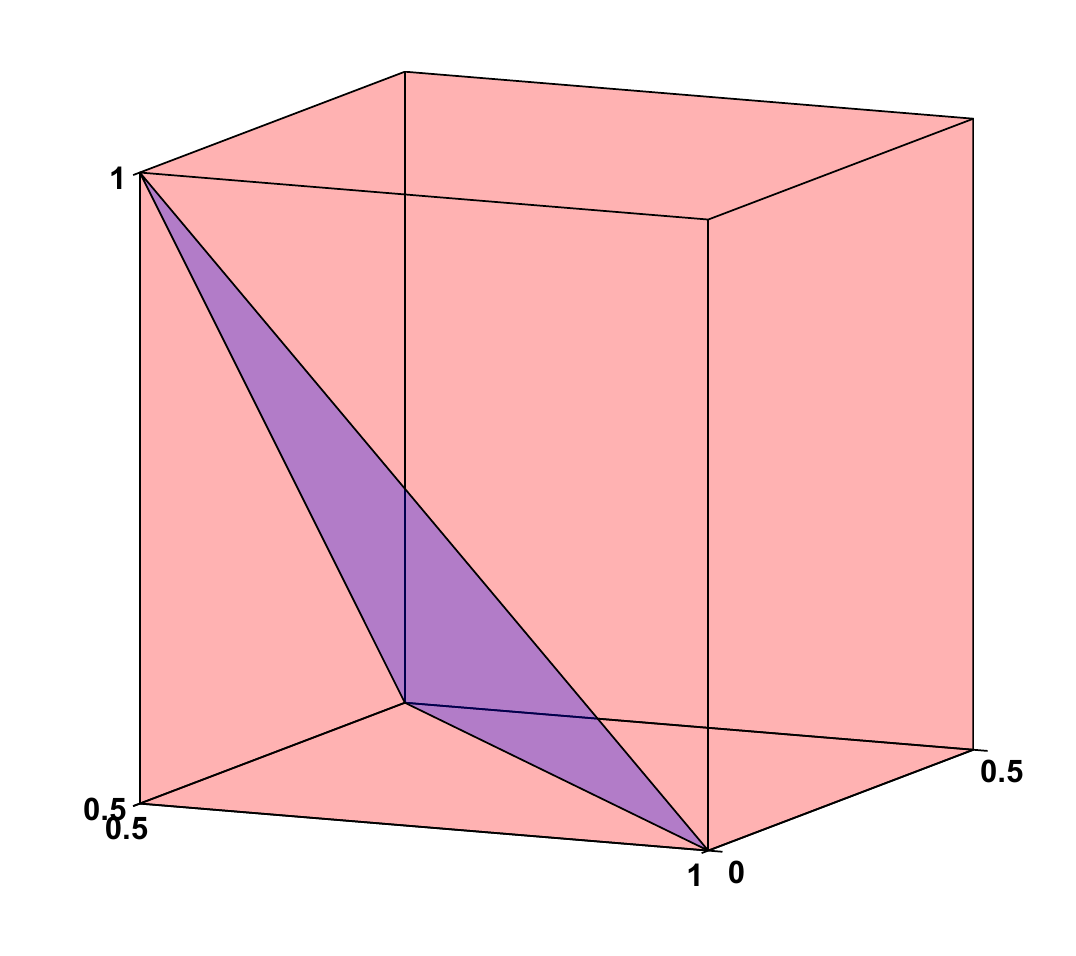}
                 \caption{Cube $8$, where $0 < q < 1/2$ and $1/2 < p,r < 1$.}
                 \end{subfigure}
                \caption{Cubes $7$ and $8$.} \label{FigSing78}
                \end{figure}

                          
As we mentioned, the map $G_{\varepsilon,3}$ is affine on each domain of continuity. The linear part is $2(1-\varepsilon)$ times the identity on each domain, and the affine part is a vector $(c_1,c_2,c_3)$ of integers, multiplied by $\varepsilon/2$. More precisely, the image of a point $(p,q,r)$ is
 \begin{align}
\left(2(1-\varepsilon)p+c_1(p,q,r) \frac{\varepsilon}{2},\text{ }2(1-\varepsilon)q+c_2(p,q,r) \frac{\varepsilon}{2},\text{ }2(1-\varepsilon)r+c_3(p,q,r) \frac{\varepsilon}{2}\right) \quad \mod 1, \label{Gc}
\end{align}
where $c_1,c_2,c_3$ take integer values depending only on the continuity domain that contains $(p,q,r)$. The exact values of $c_1,c_2$ and $c_3$ on each domain can be found in Table \ref{Tabc}.

\begin{table}[h!!]
\centering
  \begin{tabular}{| c || c | c | c | c | c || c | c || c | c | c || c | c | c |}
  	\hline
    & $\mathbf{1a}$ & $\mathbf{1b}$ & $\mathbf{1c}$ & $\mathbf{1d}$ & $\mathbf{1e}$ & $\mathbf{2a}$ & $\mathbf{2b}$ & $\mathbf{3a}$ & $\mathbf{3b}$ & $\mathbf{3c}$ & $\mathbf{4a}$ & $\mathbf{4b}$ & $\mathbf{4c}$  \\ \hline \hline
    $\mathbf{c_1}$&0&2&0&1&1&0&1&4&3&2&4&3&4 \\ \hline
    $\mathbf{c_2}$&0&1&1&0&2&4&4&4&3&3&0&1&1 \\ \hline
    $\mathbf{c_3}$&0&0&2&1&1&0&1&0&1&0&0&1&2 \\ \hline \hline
      & $\mathbf{5a}$ & $\mathbf{5b}$ & $\mathbf{5c}$& $\mathbf{6a}$ & $\mathbf{6b}$ & $\mathbf{6c}$ &$\mathbf{7a}$ & $\mathbf{7b}$ & $\mathbf{7c}$ & $\mathbf{7d}$ & $\mathbf{7e}$ & $\mathbf{8a}$ & $\mathbf{8b}$ \\ \hline \hline
      $\mathbf{c_1}$&0&1&2&0&1&0&4&2&4&3&3&4&3 \\ \hline
      $\mathbf{c_2}$&0&1&1&4&3&3&4&3&3&4&2&0&0 \\ \hline
      $\mathbf{c_3}$&4&3&4&4&3&2&4&4&2&3&3&4&3 \\ \hline    
    \end{tabular}
    \caption{Values of $c_1,c_2$ and $c_3$ in Formula \ref{Gc}, for each continuity domain.} \label{Tabc}
  \end{table}   

\section*{Appendix B: proof of Proposition 1}

The proof has three parts: we first comment on the symmetry subgroup leaving $\mathcal{A}$ invariant, then prove the dynamical invariance of $\mathcal{A}$ if the condition for $\varepsilon$ is met. Lastly, we prove that $\mathcal{A}$ is asymmetric with respect to $S_0,S_1,S_2,S_5$ and $S_6$.
  
 \textbf{\underline{Symmetries of $\mathbf{\mathcal{A}}$.}}
 We first show that $\mathcal{A}$ is symmetric with respect to $S_3$ and $S_4$. We will omit the details of the calculations, since they are straightforward and do not have any interesting details.
 
 Note that $S_3$ and $S_4$ commute. This implies that $S_4S_3(P_1)=S_3S_4(P_1)$, and by simple calculations we get that $S_4(P_2)=S_3(P_2)$. By consequence,
 \begin{align*}
 S_3S_4(P_2)&=S_4S_3(P_2)=P_2, \\
 S_3S_4S_3(P_1)&=S_4(P_1), \\
 S_4S_3S_4(P_1)&=S_3(P_1). 
 \end{align*}
(Here equalities should be understood as equalities of sets). So we can see that $S_3(\mathcal{A})= \mathcal{A}$ and $S_4(\mathcal{A})= \mathcal{A}$.

\textbf{\underline{Dynamical invariance of $\mathbf{\mathcal{A}}$.}} We start with an important observation. Suppose $S \in \langle S_3,S_4 \rangle$ and $G_{\varepsilon,3}(P_i) \subset \mathcal{A}$, $i=1$ or $2$. Then by the definition of symmetry, $G_{\varepsilon,3}(SP_i)=S G_{\varepsilon,3}(P_i) \subset S(\mathcal{A}) = \mathcal{A}$, since $\mathcal{A}$ is invariant under the symmetry subgroup $\langle S_3,S_4 \rangle$ (as a set). So it is clear that for the dynamical invariance property, we only have to check that $G_{\varepsilon,3}(P_i) \subset \mathcal{A}$ for $i=1,2$ once $\varepsilon^* \leq \varepsilon$.
 
 \underline{$G_{\varepsilon,3}(P_1) \subseteq \mathcal{A}$.} We start by showing that the image of the first polyhedron is a subset of $\mathcal{A}$, once the condition for $\varepsilon$ is met. If $ 1-\frac{\sqrt{2}}{2} \leq \varepsilon $, $P_1$ intersects two continuity domains, $1b$ and $4a$, in the following polyhedra:
 
\begin{center}
     \begin{tabular}{|p{2cm}|p{5cm}|p{5cm}|}
       \hline 
       & $\mathbf{P_1 \cap 1b}$ & $\mathbf{P_1 \cap 4a}$  \\ \hline
       $\mathbf{p}$&$p < 1/2$&$p > 1/2$ \\ \hline 
       $\mathbf{q}$&$q > \varepsilon/2$&$q > \varepsilon/2$\\ \hline 
       $\mathbf{r}$&$r > 0$&$r > 0$ \\ \hline
       $\mathbf{p+q}$&$p+q > 1-p^*$& \\ \hline
       $\mathbf{q+r}$&$q+r < p^*$& \\ \hline  
       $\mathbf{p+q+r}$&$p+q+r < 1-\varepsilon/2$&$p+q+r < 1-\varepsilon/2$ \\ \hline 
     \end{tabular}
 \end{center}
   
By applying the dynamics $G_{\varepsilon,3}$ to these sets, we get the following images:
 
\begin{center}
     \begin{tabular}{|p{2cm}|p{5cm}|p{5cm}|}
       \hline 
       & $\mathbf{G_{\varepsilon,3}(P_1 \cap 1b)}$ & $\mathbf{G_{\varepsilon,3}(P_1 \cap 4a)}$  \\ \hline
       $\mathbf{p}$&$p < 1$&$p > \varepsilon $ \\ \hline 
       $\mathbf{q}$&$q > L_{\varepsilon}(\varepsilon/2)$&$q > L_{\varepsilon}(\varepsilon/2)-\varepsilon/2$\\ \hline 
       $\mathbf{r}$&$r > 0$&$r > 0$ \\ \hline
       $\mathbf{p+q}$&$p+q > 1+p^*$& \\ \hline
       $\mathbf{q+r}$&$q+r < 1-p^*$& \\ \hline  
       $\mathbf{p+q+r}$&$p+q+r < 2-L_{\varepsilon}(\varepsilon/2)$&$p+q+r < 1-(L_{\varepsilon}(\varepsilon/2)-\varepsilon/2)$ \\ \hline 
     \end{tabular}
 \end{center}

We immediately see that $G_{\varepsilon,3}(P_1 \cap 1b) \subset P_2 \subset \mathcal{A}$.

$G_{\varepsilon,3}(P_1 \cap 4a) \subset P_1$ holds, if the conditions $p+q > 1-p^*$ and $q+r < p^*$ hold for the points $(p,q,r) \in G_{\varepsilon,3}(P_1 \cap 4a)$. 

To see when exactly (for what value of $\varepsilon$) these conditions hold, one should calculate the infimum of $p+q$ (and respectively the supremum $q+r$) on the polyhedron $G_{\varepsilon,3}(P_1 \cap 4a)$. In other words, one should minimze the objective function $c_1=p+q$ (maximize $c_2=q+r$). This is a linear optimization task. To solve it, one could apply the simplex algorithm for example, but in this simple case we can also solve it by a trial and error method: we know that the optimal value will be attained at an extremal point (a vertice) of the polyhedron (which is a tetrahedron in fact), so we only have to compare the value of $p+q$ (or $q+r$, respectively) in four places.

This way we obtain that the infimum of $p+q$ on $G_{\varepsilon,3}(P_1 \cap 4a)$ is $L_{\varepsilon}(\varepsilon/2)+\varepsilon/2$, so for $G_{\varepsilon,3}(P_1 \cap 4a) \subset P_1$ 
to hold we must have
\begin{align}
L_{\varepsilon}(\varepsilon/2)+\varepsilon/2 & \geq 1-p^*, \nonumber \\
p^* & \geq (1-\varepsilon)^2, \nonumber \\
\varepsilon & \geq \varepsilon^* \label{epsstar1}.
\end{align}
\begin{figure}
 \centering
 \begin{subfigure}[b]{0.3\textwidth}
 \includegraphics[scale=0.42]{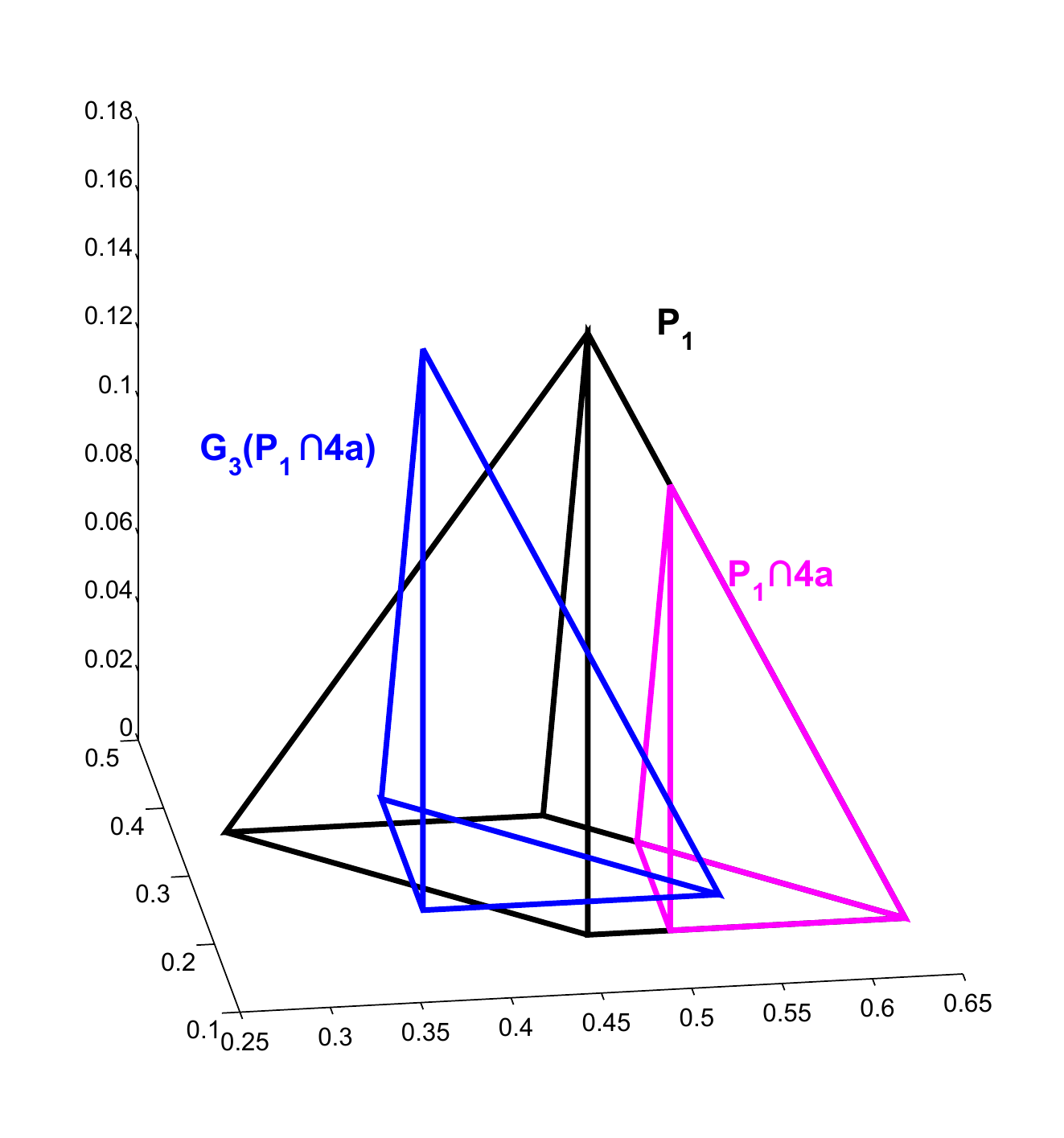}
 \caption{$\varepsilon=0.37$}
 \end{subfigure}
 \hspace{2cm}
 \begin{subfigure}[b]{0.3\textwidth}
  \includegraphics[scale=0.45]{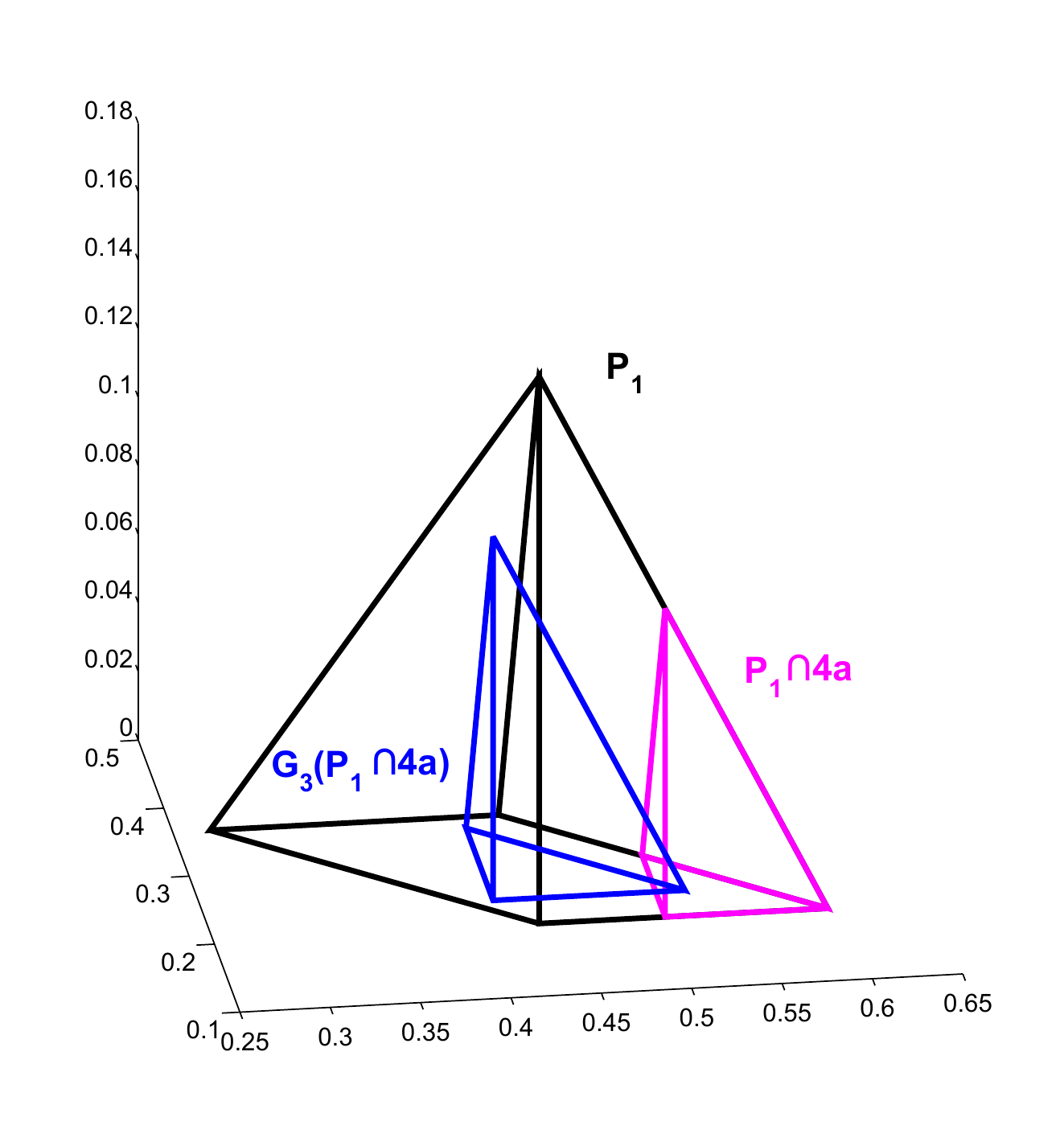}
  \caption{$\varepsilon=0.41$}
  \end{subfigure}
 \caption{The image of $P_1 \cap 4a$ for $\varepsilon=0.37 < \varepsilon^*$ and $\varepsilon=0.41 > \varepsilon^*$.} \label{FigErgB1}
 \end{figure}
Similarly, the supremum of $q+r$ on $G_{\varepsilon,3}(P_1 \cap 4a)$ is $1-L_{\varepsilon}(\varepsilon/2)-\varepsilon/2$, so for $G_{\varepsilon,3}(P_1 \cap 4a) \subset P_1$ to hold we must have
\begin{align}
1-L_{\varepsilon}(\varepsilon/2) - \varepsilon/2 & \leq p^*, \nonumber \\
1-p^* & \leq L_{\varepsilon}(\varepsilon/2)+\varepsilon/2 , \nonumber\\
(1-\varepsilon)^2 & \leq p^*, \nonumber \\
\varepsilon^* & \leq \varepsilon \label{epsstar2}. 
\end{align}
For an illustration, see Figure \ref{FigErgB1}. 

\underline{$G_{\varepsilon,3}(P_2) \subseteq \mathcal{A}$.} Similar calculations are needed to show this for appropriate values of $\varepsilon$. The set $P_2$ intersects six continuity domains, $3a$, $3b$, $3c$, $4a$, $4b$ and $4c$. Computing the images of these intersections, one can conclude that 
\begin{align*}
G_{\varepsilon,3}(P_2 \cap 3b) &\subset S_3(P_1) \\
G_{\varepsilon,3}(P_2 \cap 3c) &\subset P_1 \\
G_{\varepsilon,3}(P_2 \cap 4b) &\subset S_4(P_1) \\
G_{\varepsilon,3}(P_2 \cap 3c) &\subset S_3S_4(P_1)
\end{align*}  
always hold, and
\begin{align*}
G_{\varepsilon,3}(P_2 \cap 3a) &\subset P_2 \\
G_{\varepsilon,3}(P_2 \cap 4a) &\subset P_2
\end{align*}
hold if and only if $\varepsilon \geq \frac{1}{8}(7-\sqrt{17}) \approx 0.359$. Now we see that $\frac{1}{8}(7-\sqrt{17}) < \varepsilon^*$, so $G_{\varepsilon,3}(P_2) \subseteq \mathcal{A}$ holds if $\varepsilon \geq \varepsilon^*$.

 \textbf{\underline{Asymmetries of $\mathbf{\mathcal{A}}$.}} 

\begin{figure}[h!]
      \centering
      \includegraphics[scale=0.6]{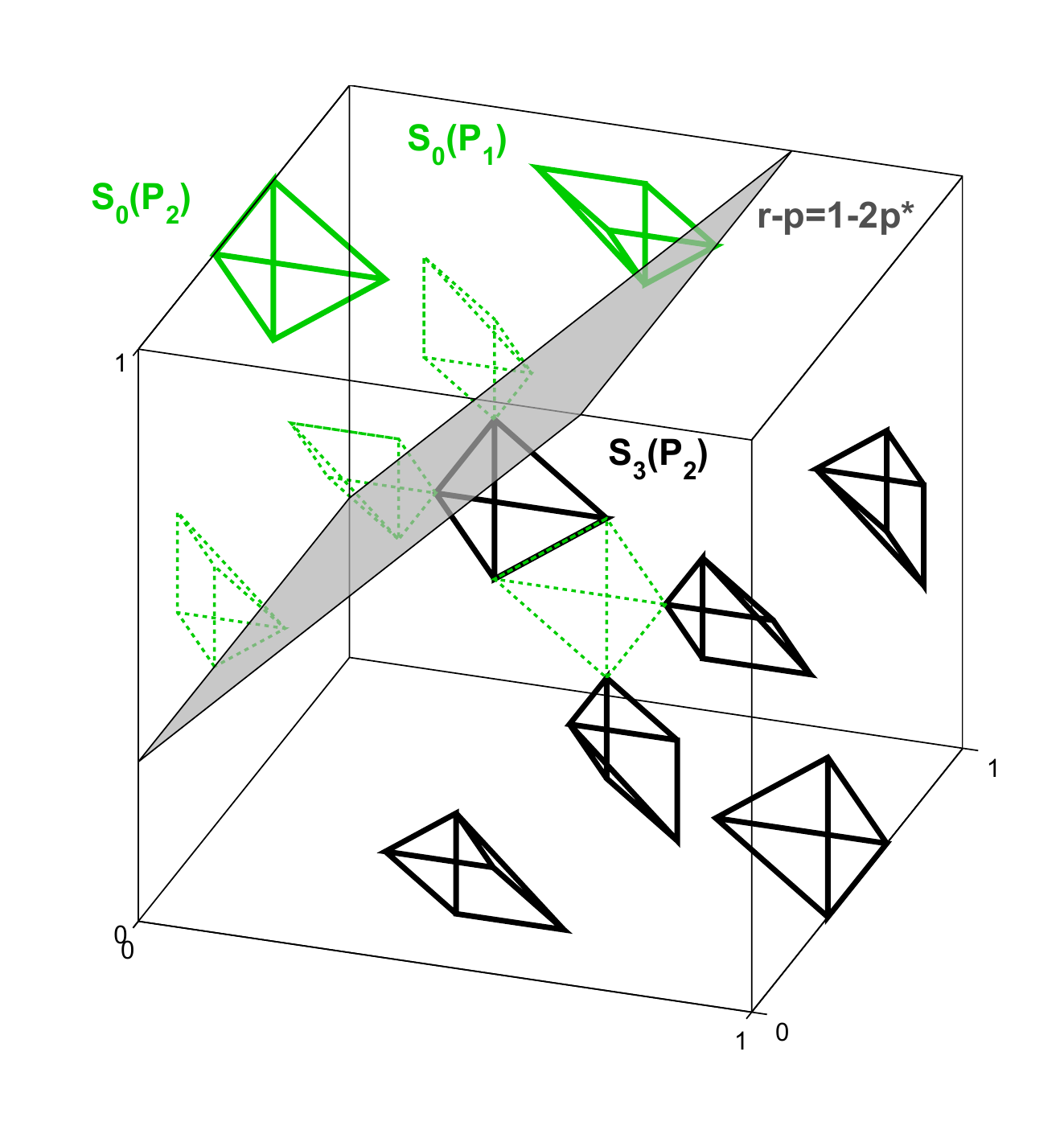}
      \caption{The sets $\mathcal{A}$ (black) and $S_0(\mathcal{A})$ (green). The plane $r-p=1-2p^*$ separates $S_0(P_1)$, $S_0(P_2)$ and $\mathcal{A}$.} \label{FigS0}
      \end{figure}
      
We first prove the asymmetry with respect to $S_0$, derived from the inversion symmetry of $F_{\varepsilon,3}$, then prove asymmetry with respect to $S_1,S_2,S_5$ and $S_6$ derived from the permutation symmetries of $F_{\varepsilon,3}$.

\underline{Asymmetry with respect to $S_0$.} Notice that it is enough to check that both $S_0(P_1)$ and  $S_0(P_2)$ are disjoint from $\mathcal{A}$. This is enough indeed, since if $S_0(P_i) \subset \mathcal{A}^C$, $i=1,2$, then $S_0S(P_i)=SS_0(P_i) \subset S\mathcal{A}^C = \mathcal{A}^C$, if $S \subset \langle S_3,S_4 \rangle$.

In the table below, we describe $S_0(P_1), S_0(P_2)$ and $S_3(P_2) \subset \mathcal{A}$.

 \begin{center}
        \begin{tabular}{|p{2cm}|p{3cm}|p{3cm}|p{3cm}|}
          \hline 
          & $\mathbf{S_0(P_1)}$ & $\mathbf{S_0(P_2)}$& $\mathbf{S_3(P_2)}$    \\ \hline
          $\mathbf{p}$&&$p > 0$& $p > p^*$ \\ \hline 
          $\mathbf{q}$&$q < 1-\varepsilon/2$&& \\ \hline 
          $\mathbf{r}$&$r < 1$&$r < 1$& $r < 1-p^*$ \\ \hline
          $\mathbf{p+q}$&$p+q < 1+p^*$&$p+q < 1-p^*$& $p+q < 1$ \\ \hline
          $\mathbf{q+r}$&$q+r > 2-p^*$&$q+r > 1+p^*$& $q+r > 1$ \\ \hline  
          $\mathbf{p+q+r}$&$p+q+r > 2+\varepsilon/2$&& \\ \hline 
        \end{tabular}
    \end{center}    
 
We are going to show that the plane $r-p=1-2p^*$ separates both $S_0(P_1)$ and $S_0(P_2)$ from $\mathcal{A}$. We see that the infimum of $r-p$ on $S_0(P_1)$ is exactly $1-2p^*$, and it is $2p^* > 1-2p^*$ on $S_0(P_2)$ (it is easy to see that $p^* \geq 1/3$ for every value of $\varepsilon$). So these sets are above the plane. On the other hand, the supremum of $r-p$ on $\mathcal{A}$ is $1-2p^*$ (attained on the closure of $S_3(P_2)$), hence this plane separates $\mathcal{A}$ from $S_0(P_1)$ and $S_0(P_2)$. For an illustration see Figure \ref{FigS0}.
 
\underline{Asymmetry with respect to $S_1,S_2,S_5$ and $S_6$.} We remind the reader that $S_5=S_3S_1S_3$, $S_2=S_4S_1S_4$ and $S_6=S_3S_2S_3$. Notice that it is enough to prove that $S_1(\mathcal{A}) \subset \mathcal{A}^C$. From this it follows that 
 \begin{align*}
 S_5(\mathcal{A})&=S_3S_1(\mathcal{A}) \subset S_3\mathcal{A}^C = \mathcal{A}^C, \\
 S_2(\mathcal{A})&=S_4S_1(\mathcal{A}) \subset S_4\mathcal{A}^C = \mathcal{A}^C, \\
 S_6(\mathcal{A})&=S_3S_4S_1(\mathcal{A}) \subset S_3S_4\mathcal{A}^C =\mathcal{A}^C.
 \end{align*}
 We now describe the six polyhedra of $S_1(\mathcal{A})$.
 \begin{center}
          \begin{tabular}{|p{2cm}|p{3cm}|p{3cm}| p{3cm}|}
            \hline 
            & $\mathbf{S_1(P_1)}$ & $\mathbf{S_1(P_2)}$& $\mathbf{S_1S_3(P_1)}$    \\ \hline
            $\mathbf{p}$&&$p > 0$& $p < p^*$ \\ \hline 
            $\mathbf{q}$&$q > 1-p^*$&$q > p^*$& $q < 1$ \\ \hline 
            $\mathbf{r}$&$r > 0$&$r > 0$&$r < p^*$ \\ \hline
            $\mathbf{p+q}$&$p+q > 1+\varepsilon/2$&& $p+q > 1-\varepsilon/2$ \\ \hline
            $\mathbf{q+r}$&$q+r < 1-\varepsilon/2$&& $q+r > \varepsilon/2$ \\ \hline  
            $\mathbf{p+q+r}$&$p+q+r < 1+p^*$&$p+q+r < 1-p^*$& \\ \hline 
          \end{tabular}
          \vspace{0.5cm} \\ 
          \begin{tabular}{|p{2cm}|p{3cm}|p{3cm}| p{3cm}|}
                      \hline 
                      & $\mathbf{S_1S_4(P_1)}$ & $\mathbf{S_1S_3(P_2)}$& $\mathbf{S_1S_3S_4(P_1)}$    \\ \hline
                      $\mathbf{p}$&$p < p^*$&$p < 1-p^*$& $p > 0$ \\ \hline 
                      $\mathbf{q}$&&$q < 1$&$q > 1-p^*$ \\ \hline 
                      $\mathbf{r}$&$r < p^*$&$r<1-p^*$& \\ \hline
                      $\mathbf{p+q}$&$p+q < 1-\varepsilon/2$&&$p+q < 1-\varepsilon/2$ \\ \hline
                      $\mathbf{q+r}$&$q+r < 1-\varepsilon/2$&&$q+r > 1-\varepsilon/2$ \\ \hline  
                      $\mathbf{p+q+r}$&$p+q+r > 1$&$p+q+r > 2$&$p+q+r < 1+p^*$ \\ \hline 
                    \end{tabular}
      \end{center}
      
%
%
Similar calculations as in the previous point lead to the observation that 
\begin{table}[h!!]
\centering
\begin{tabular}{llll}
the plane & $p+q+r=1-p^*$ & separates & $S_1(P_2)$ and $\mathcal{A},$ \\
the plane & $q+r=1+p^*$ & separates & $S_1S_3(P_2)$ and $\mathcal{A}.$
\end{tabular}
\end{table}

\begin{figure}[h!]
        \centering
        
        \begin{subfigure}[b]{0.4\textwidth}
        \includegraphics[scale=0.55]{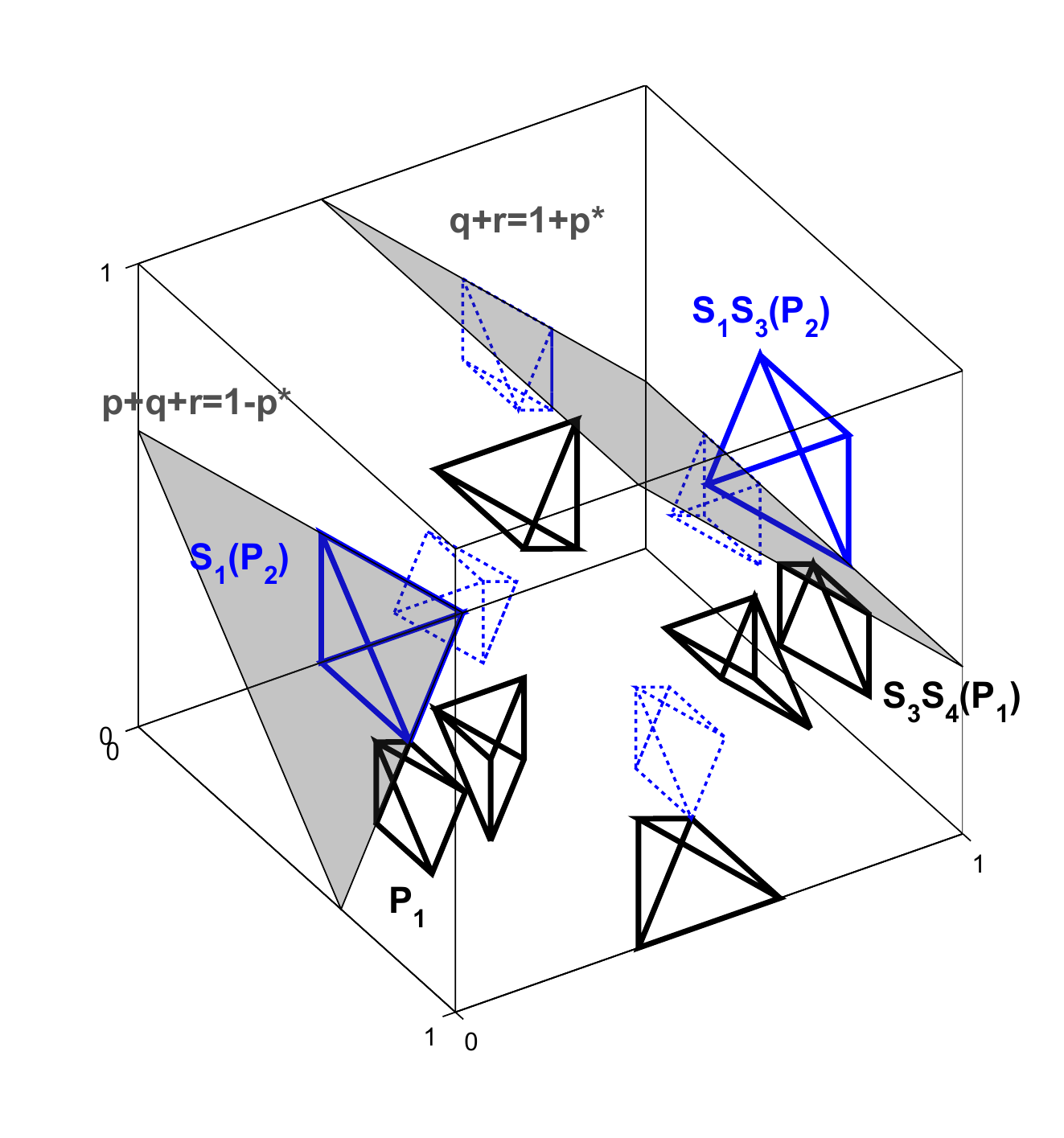}
        \caption{The plane $p+q+r=1-p^*$ separates $S_1(P_2)$ from $\mathcal{A}$ and the plane $q+r=1+p^*$ separates $S_1S_3(P_2)$ from $\mathcal{A}$.} \label{FigS1a}
        \end{subfigure}
                 \hspace{2cm}
        \centering
        \begin{subfigure}[b]{0.4\textwidth}
         \includegraphics[scale=0.5]{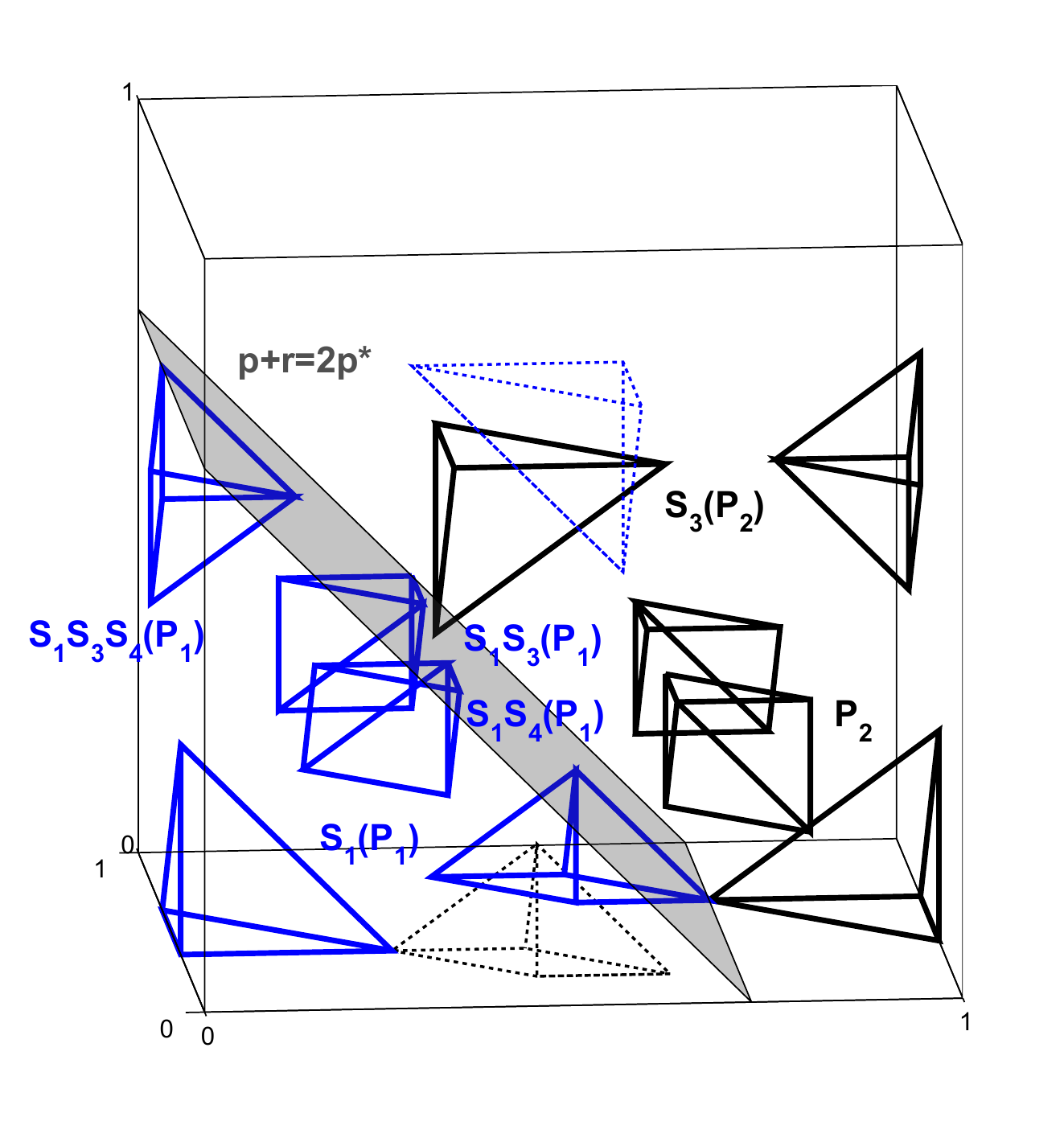}
         \caption{The plane $p+r=2p^*$ separates $S_1(P_1), S_1S_3(P_1), S_1S_4(P_1)$ and $S_1S_3S_4(P_1)$ from $\mathcal{A} \backslash P_1$ .} \label{FigS1b}
         \end{subfigure}
         \hspace{2cm}
         \begin{subfigure}[b]{0.4\textwidth}
                 \includegraphics[scale=0.55]{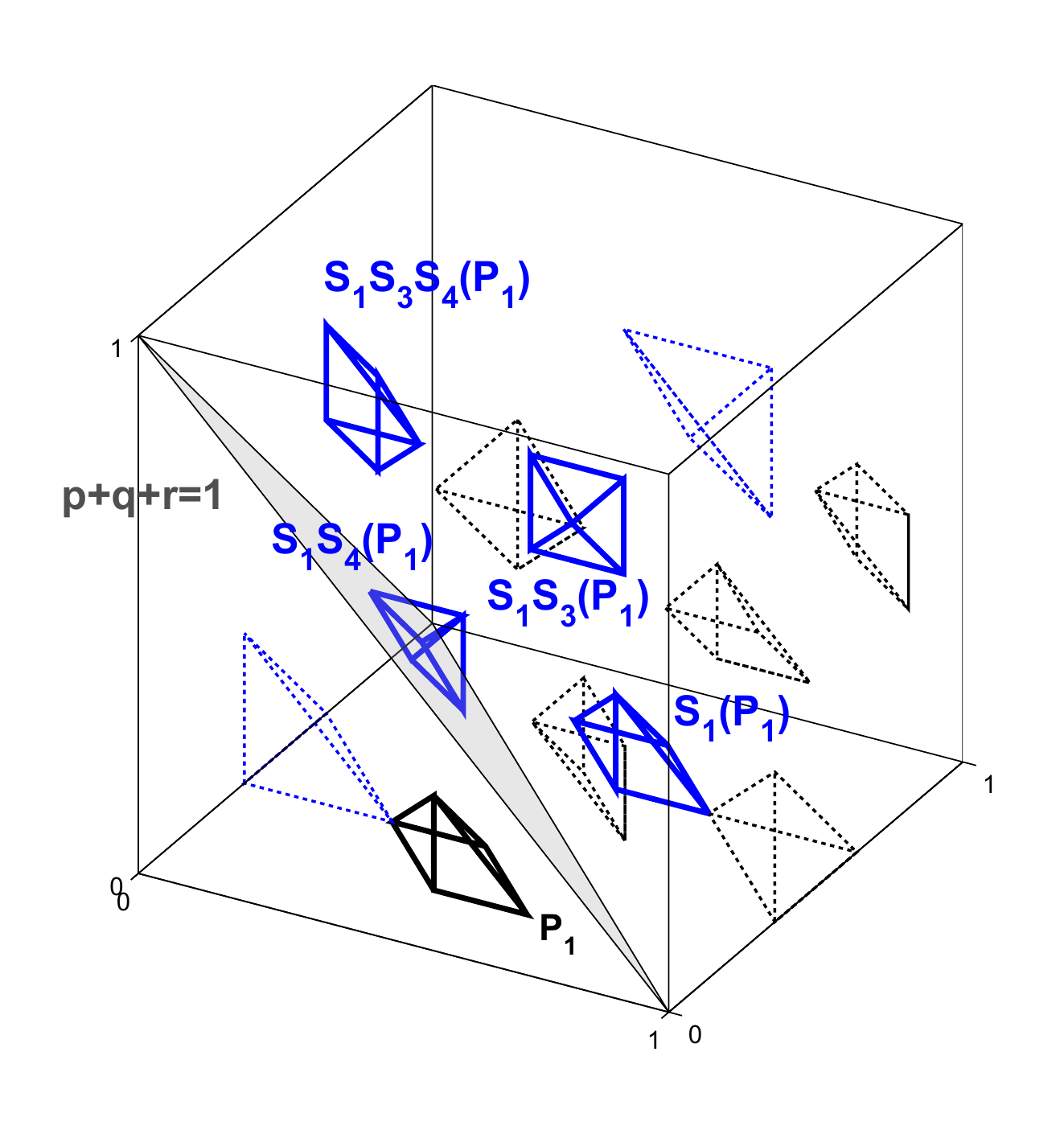}
                 \caption{The plane $p+q+r=1$ separates $S_1(P_1), S_1S_3(P_1), S_1S_4(P_1)$ and $S_1S_3S_4(P_1)$ from $P_1$ .} \label{FigS1c}
                 \end{subfigure}
        \caption{The sets $\mathcal{A}$ (black) and $S_1(\mathcal{A})$ (blue). Different angles are plotted for better visibility.} \label{FigS1}
        \end{figure}

For an illustration, see \Cref{FigS1a}.
 
Furthermore, 
\begin{table}[h!!]
\centering
\begin{tabular}{llll}
the plane & $p+r=2p^*$ & separates & $S_1(P_1) \cup S_1S_3(P_1) \cup S_1S_4(P_1) \cup S_1S_3S_4(P_1)$ and $\mathcal{A} \backslash P_1,$ \\
the plane & $p+q+r=1$ & separates & $S_1(P_1) \cup S_1S_3(P_1) \cup S_1S_4(P_1) \cup S_1S_3S_4(P_1)$ and $P_1.$
\end{tabular}
\end{table}

For an illustration, see \Cref{FigS1b} and  \Cref{FigS1c}.

\section*{Appendix C: Proof of Proposition 2}

To check the symmetry properties, many straightforward calculations are needed (it suffices to check that the set $\mathcal{S}$ is symmetric with respect to the minimal generating symmetries $S_0, S_1, S_2$ and $ S_3$). We omit these calculations, since they do not bear any interesting details.

We now prove that $G_{\varepsilon,3}(\mathcal{S}) \subseteq \mathcal{S}$ if $\varepsilon \geq 1-\frac{\sqrt{2}}{2}$. To obtain this, it is enough to check that $G_{\varepsilon,3}(P_0) \subseteq \mathcal{S}$ by similar arguments used to show that the invariance of $\mathcal{A}$ was a consequence of $G_{\varepsilon,3}(P_i) \subseteq \mathcal{A}$, $i=1,2$.


First of all, we are going to describe the 12 polyhedra of $\mathcal{S}$ in a table below. We will give descriptions with inequalities which are always correct in $\mathbb{R}^3$ and correct in $\mathbb{T}^3$, if $1-\frac{\sqrt{2}}{2} \leq \varepsilon$.  
  
\begin{center}
    \begin{tabular}{|p{2cm}|p{7cm}|p{7cm}|}
      \hline 
      & $\mathbf{P_0}$ & $\mathbf{S_0(P_0)}$  \\ \hline
      $\mathbf{p}$ &$L_{\varepsilon}(\varepsilon/2) < p < L_{\varepsilon}(1-\varepsilon/2)$ & $L(\varepsilon/2) < p < L_{\varepsilon}(1-\varepsilon/2)$ \\ \hline
      $\mathbf{q}$ &$\varepsilon/2 < q < L_{\varepsilon}^2(1-\varepsilon/2)$& $L_{\varepsilon}^2(\varepsilon/2) < q < 1-\varepsilon/2$ \\ \hline
      $\mathbf{r}$ & $L_{\varepsilon}(\varepsilon/2) < r < L_{\varepsilon}(1-\varepsilon/2)$ & $L_{\varepsilon}(\varepsilon/2) < r < L_{\varepsilon}(1-\varepsilon/2)$ \\ \hline
      $\mathbf{p+q}$ && \\ \hline
      $\mathbf{q+r}$ && \\ \hline
      $\mathbf{p+q+r}$ & $1+\varepsilon/2 < p+q+r < 1+L_{\varepsilon}^2(1-\varepsilon/2)$ & $1+L_{\varepsilon}^2(\varepsilon/2) <p+q+r <2-\varepsilon/2$ \\ \hline
    \end{tabular}
\end{center}

\begin{center}
    \begin{tabular}{|p{2cm}|p{7cm}|p{7cm}|}
      \hline 
      & $\mathbf{S_1(P_0)}$ & $\mathbf{S_2(P_0)}$  \\ \hline
      $\mathbf{p}$ &$L_{\varepsilon}(\varepsilon/2) < p < L_{\varepsilon}(1-\varepsilon/2)$ & $L_{\varepsilon}^2(\varepsilon/2) < p < 1-\varepsilon/2$ \\ \hline
      $\mathbf{q}$ && $L_{\varepsilon}(\varepsilon/2) < q < L_{\varepsilon}(1-\varepsilon/2)$ \\ \hline
      $\mathbf{r}$ & $L_{\varepsilon}(\varepsilon/2) < r < L_{\varepsilon}(1-\varepsilon/2)$ & $\varepsilon/2 < r < L_{\varepsilon}^2(1-\varepsilon/2)$ \\ \hline
      $\mathbf{p+q}$ & $1+\varepsilon/2 < p+q < 1+L_{\varepsilon}^2(1-\varepsilon/2)$& \\ \hline
      $\mathbf{q+r}$ & $1+\varepsilon/2 < q+r < 1+L_{\varepsilon}^2(1-\varepsilon/2)$ & \\ \hline
      $\mathbf{p+q+r}$ & & $1+L_{\varepsilon}(\varepsilon/2) <p+q+r <1+L_{\varepsilon}(1-\varepsilon/2)$ \\ \hline
    \end{tabular}
\end{center}

\begin{center}
    \begin{tabular}{|p{2cm}|p{7cm}|p{7cm}|}
      \hline 
      & $\mathbf{S_3(P_0)}$ & $\mathbf{S_4(P_0)}$  \\ \hline
      $\mathbf{p}$ & &  \\ \hline
      $\mathbf{q}$ & $\varepsilon/2 < q < L_{\varepsilon}^2(1-\varepsilon/2)$ & $\varepsilon/2 < q < L_{\varepsilon}^2(1-\varepsilon/2)$ \\ \hline
      $\mathbf{r}$ &  &  \\ \hline
      $\mathbf{p+q}$ & $L_{\varepsilon}(\varepsilon/2) < p+q < L_{\varepsilon}(1-\varepsilon/2)$& $1+L_{\varepsilon}(\varepsilon/2) < p+q < 1+L_{\varepsilon}(1-\varepsilon/2)$ \\ \hline
      $\mathbf{q+r}$ & $L_{\varepsilon}(\varepsilon/2) < q+r < L_{\varepsilon}(1-\varepsilon/2)$ & $1+L_{\varepsilon}(\varepsilon/2) < q+r < 1+L_{\varepsilon}(1-\varepsilon/2)$ \\ \hline
      $\mathbf{p+q+r}$ & $L_{\varepsilon}^2(\varepsilon/2) < p+q+r < 1-\varepsilon/2$ & $2+\varepsilon/2 <p+q+r <2+L_{\varepsilon}^2(1-\varepsilon/2)$ \\ \hline
    \end{tabular}
\end{center}

\begin{center}
    \begin{tabular}{|p{2cm}|p{7cm}|p{7cm}|}
      \hline 
      & $\mathbf{S_5(P_0)}$ & $\mathbf{S_0S_1(P_0)}$  \\ \hline
      $\mathbf{p}$ & $\varepsilon/2 < p < L_{\varepsilon}^2(1-\varepsilon/2)$& $L_{\varepsilon}(\varepsilon/2) < p < L_{\varepsilon}(1-\varepsilon/2)$  \\ \hline
      $\mathbf{q}$ & $L_{\varepsilon}(\varepsilon/2) < q < L_{\varepsilon}(1-\varepsilon/2)$ &  \\ \hline
      $\mathbf{r}$ & $L_{\varepsilon}^2(\varepsilon/2) < r < 1-\varepsilon/2$ & $L_{\varepsilon}(\varepsilon/2) < r < L_{\varepsilon}(1-\varepsilon/2)$ \\ \hline
      $\mathbf{p+q}$ & & $L_{\varepsilon}^2(\varepsilon/2) < p+q < 1-\varepsilon/2$ \\ \hline
      $\mathbf{q+r}$ &  & $L_{\varepsilon}^2(\varepsilon/2) < q+r < 1-\varepsilon/2$ \\ \hline
      $\mathbf{p+q+r}$ & $1+L_{\varepsilon}(\varepsilon/2) < p+q+r < 1+L_{\varepsilon}(1-\varepsilon/2)$ &  \\ \hline
    \end{tabular}
\end{center}

\begin{center}
    \begin{tabular}{|p{2cm}|p{7cm}|p{7cm}|}
      \hline 
      & $\mathbf{S_2S_1(P_0)}$ & $\mathbf{S_3S_1(P_0)}$  \\ \hline
      $\mathbf{p}$ & & $L_{\varepsilon}^2(\varepsilon/2) < p < 1-\varepsilon/2$  \\ \hline
      $\mathbf{q}$ & $L_{\varepsilon}(\varepsilon/2) < q < L_{\varepsilon}(1-\varepsilon/2)$ &  \\ \hline
      $\mathbf{r}$ &  & $L_{\varepsilon}^2(\varepsilon/2) < r < 1-\varepsilon/2$  \\ \hline
      $\mathbf{p+q}$ & $L_{\varepsilon}^2(\varepsilon/2) < p+q < 1-\varepsilon/2$ & $1+L_{\varepsilon}(\varepsilon/2) < p+q < 1+L_{\varepsilon}(1-\varepsilon/2)$ \\ \hline
      $\mathbf{q+r}$ & $1+\varepsilon/2 < q+r < 1+L_{\varepsilon}^2(1-\varepsilon/2)$ & $1+L_{\varepsilon}(\varepsilon/2) < q+r < 1+L_{\varepsilon}(1-\varepsilon/2)$ \\ \hline
      $\mathbf{p+q+r}$ & $1+L_{\varepsilon}(\varepsilon/2) < p+q+r < 1+L_{\varepsilon}(1-\varepsilon/2)$ &  \\ \hline
    \end{tabular}
\end{center}

\begin{center}
    \begin{tabular}{|p{2cm}|p{7cm}|p{7cm}|}
      \hline 
      & $\mathbf{S_4S_1(P_0)}$ & $\mathbf{S_5S_1(P_0)}$  \\ \hline
      $\mathbf{p}$ & $\varepsilon/2 < p < L_{\varepsilon}^2(1-\varepsilon/2)$ &   \\ \hline
      $\mathbf{q}$ &  & $L_{\varepsilon}(\varepsilon/2) < q < L_{\varepsilon}(1-\varepsilon/2)$  \\ \hline
      $\mathbf{r}$ & $\varepsilon/2 < r < L_{\varepsilon}^2(1-\varepsilon/2)$ &   \\ \hline
      $\mathbf{p+q}$ & $L_{\varepsilon}(\varepsilon/2) < p+q < L_{\varepsilon}(1-\varepsilon/2)$ & $1+\varepsilon/2 < p+q < 1+L_{\varepsilon}^2(1-\varepsilon/2)$ \\ \hline
      $\mathbf{q+r}$ & $L_{\varepsilon}(\varepsilon/2) < q+r < L_{\varepsilon}(1-\varepsilon/2)$ & $L_{\varepsilon}^2(\varepsilon/2) < q+r < 1-\varepsilon/2$ \\ \hline
      $\mathbf{p+q+r}$ &  & $1+L_{\varepsilon}(\varepsilon/2) < p+q+r < 1+L_{\varepsilon}(1-\varepsilon/2)$ \\ \hline
    \end{tabular}
\end{center}

Observe that the polyhedron $P_0$ intersects four domains of continuity, $1e$, $4b$, $5b$ and $8b$ in the following way:\begin{center}
    \begin{tabular}{|p{2cm}|p{7cm}|p{7cm}|}
      \hline 
      &$\mathbf{P_0 \cap 1e}$ &  $\mathbf{P_0 \cap 4b}$ \\ \hline
      $\mathbf{p}$ & $L_{\varepsilon}(\varepsilon/2) < p < 1/2$ & $1/2 < p < L_{\varepsilon}(1-\varepsilon/2)$ \\ \hline
      $\mathbf{q}$ & $\varepsilon/2 < q < L_{\varepsilon}^2(1-\varepsilon/2)$ & $\varepsilon/2 < q < L_{\varepsilon}^2(1-\varepsilon/2)$ \\ \hline 
      $\mathbf{r}$ &$L_{\varepsilon}(\varepsilon/2) < r < 1/2$ &  $L_{\varepsilon}(\varepsilon/2) < r < 1/2$ \\ \hline
      $\mathbf{p+q+r}$ & $1+\varepsilon/2 < p+q+r < 1+L_{\varepsilon}^2(1-\varepsilon/2)$ & $1+\varepsilon/2 < p+q+r < 1+L_{\varepsilon}^2(1-\varepsilon/2)$ \\ \hline  
    \end{tabular}
\end{center}
\begin{center}
    \begin{tabular}{|p{2cm}|p{7cm}|p{7cm}|}
      \hline 
      & $\mathbf{P_0 \cap 5b} $ & $\mathbf{P_0 \cap 8b}$ \\ \hline
      $\mathbf{p}$&$L_{\varepsilon}(\varepsilon/2) < p < 1/2$& $1/2 < p < L_{\varepsilon}(1-\varepsilon/2)$ \\ \hline 
      $\mathbf{q}$&$\varepsilon/2 < q < L_{\varepsilon}^2(1-\varepsilon/2)$&$\varepsilon/2 < q < L_{\varepsilon}^2(1-\varepsilon/2)$ \\ \hline 
      $\mathbf{r}$&$1/2 < r < L_{\varepsilon}(1-\varepsilon/2)$& $1/2 < r < L_{\varepsilon}(1-\varepsilon/2)$ \\ \hline  
      $\mathbf{p+q+r}$&$1+\varepsilon/2 < p+q+r < 1+L_{\varepsilon}^2(1-\varepsilon/2)$& $1+\varepsilon/2 < p+q+r < 1+L_{\varepsilon}^2(1-\varepsilon/2)$ \\ \hline 
    \end{tabular}
\end{center}

We calculate the images of these polyhedra, using that the coordinates $p,q$ and $r$ evolve according to $L_{\varepsilon}$ in domains $4b$ and $5b$ and the coordinates $p$ and $r$ evolve according to $L_{\varepsilon}$ in domains $1e$ and $8b$. The results are collected in the table below.

\begin{center}
    \begin{tabular}{|p{2cm}|p{7cm}|p{7cm}|}
      \hline 
      &$\mathbf{G_{\varepsilon,3}(P_0 \cap 1e)}$ &  $\mathbf{G_{\varepsilon,3}(P_0 \cap 4b)}$ \\ \hline
      $\mathbf{p}$ & $L_{\varepsilon}^2(\varepsilon/2) < p < 1-\varepsilon/2$ & $\varepsilon/2 < p < L_{\varepsilon}^2(1-\varepsilon/2)$ \\ \hline
      $\mathbf{q}$ & $L_{\varepsilon}(\varepsilon/2)+\varepsilon/2 < q < L_{\varepsilon}^3(1-\varepsilon/2)+\varepsilon/2$ & $L_{\varepsilon}(\varepsilon/2) < q < L_{\varepsilon}^3(1-\varepsilon/2)$ \\ \hline 
      $\mathbf{r}$ &$L_{\varepsilon}^2(\varepsilon/2) < r < 1-\varepsilon/2$ &  $L_{\varepsilon}^2(\varepsilon/2) < r < 1-\varepsilon/2$ \\ \hline
      $\mathbf{p+q+r}$ & $1+L_{\varepsilon}(\varepsilon/2)+\varepsilon/2 < p+q+r < 1+L_{\varepsilon}^3(1-\varepsilon/2)+\varepsilon/2$ & $1+L_{\varepsilon}(\varepsilon/2) < p+q+r < 1+L_{\varepsilon}^3(1-\varepsilon/2)$ \\ \hline  
    \end{tabular}
\end{center}
\begin{center}
    \begin{tabular}{|p{2cm}|p{7cm}|p{7cm}|}
      \hline 
      & $\mathbf{G_{\varepsilon,3}(P_0 \cap 5b)} $ & $\mathbf{G_{\varepsilon,3}(P_0 \cap 8b)}$ \\ \hline
      $\mathbf{p}$&$L_{\varepsilon}^2(\varepsilon/2) < p < 1-\varepsilon/2$& $\varepsilon/2 < p < L_{\varepsilon}^2(1-\varepsilon/2)$ \\ \hline 
      $\mathbf{q}$&$L_{\varepsilon}(\varepsilon/2) < q < L_{\varepsilon}^3(1-\varepsilon/2)$&$L_{\varepsilon}(\varepsilon/2)-\varepsilon/2 < q < L_{\varepsilon}^3(1-\varepsilon/2) -\varepsilon/2$ \\ \hline 
      $\mathbf{r}$&$\varepsilon/2 < r < L_{\varepsilon}^2(1-\varepsilon/2)$& $\varepsilon/2 < r < L_{\varepsilon}^2(1-\varepsilon/2)$ \\ \hline  
      $\mathbf{p+q+r}$&$1+L_{\varepsilon}(\varepsilon/2) < p+q+r < 1+L_{\varepsilon}^3(1-\varepsilon/2)$& $1+L_{\varepsilon}(\varepsilon/2)-\varepsilon/2 < p+q+r < 1+L_{\varepsilon}^3(1-\varepsilon/2)-\varepsilon/2$ \\ \hline 
    \end{tabular}
\end{center}

We immediately see that $G_{\varepsilon,3}(P_0 \cap 4b) \subseteq S_5(P_0)$ if and only if 
\begin{align}
L_{\varepsilon}^3(1-\varepsilon/2) & \leq L_{\varepsilon}(1-\varepsilon/2) \nonumber \\
1-\frac{\sqrt{2}}{2} &\leq \varepsilon. \label{sqrt2}
\end{align}
 
 \begin{figure}[h!]
  \centering
  \begin{subfigure}[b]{0.3\textwidth}
  \includegraphics[scale=0.45]{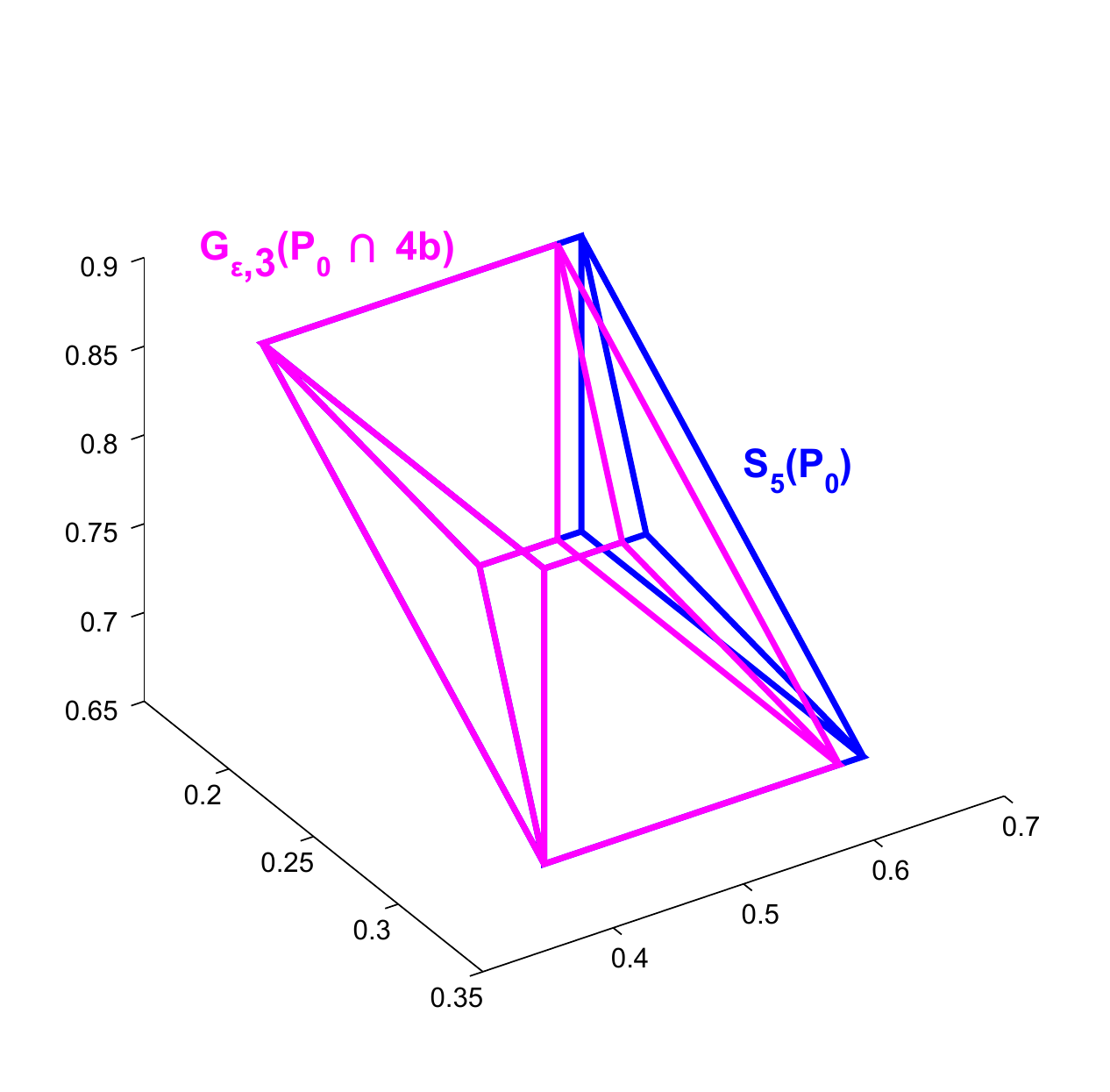}
  \caption{The image of $P_0 \cap 4b$. The case of the image of $P_0 \cap 5b$ is geometrically very similar.} \label{Fig45b}
  \end{subfigure}
  \hspace{2cm}
  \begin{subfigure}[b]{0.3\textwidth}
   \includegraphics[scale=0.45]{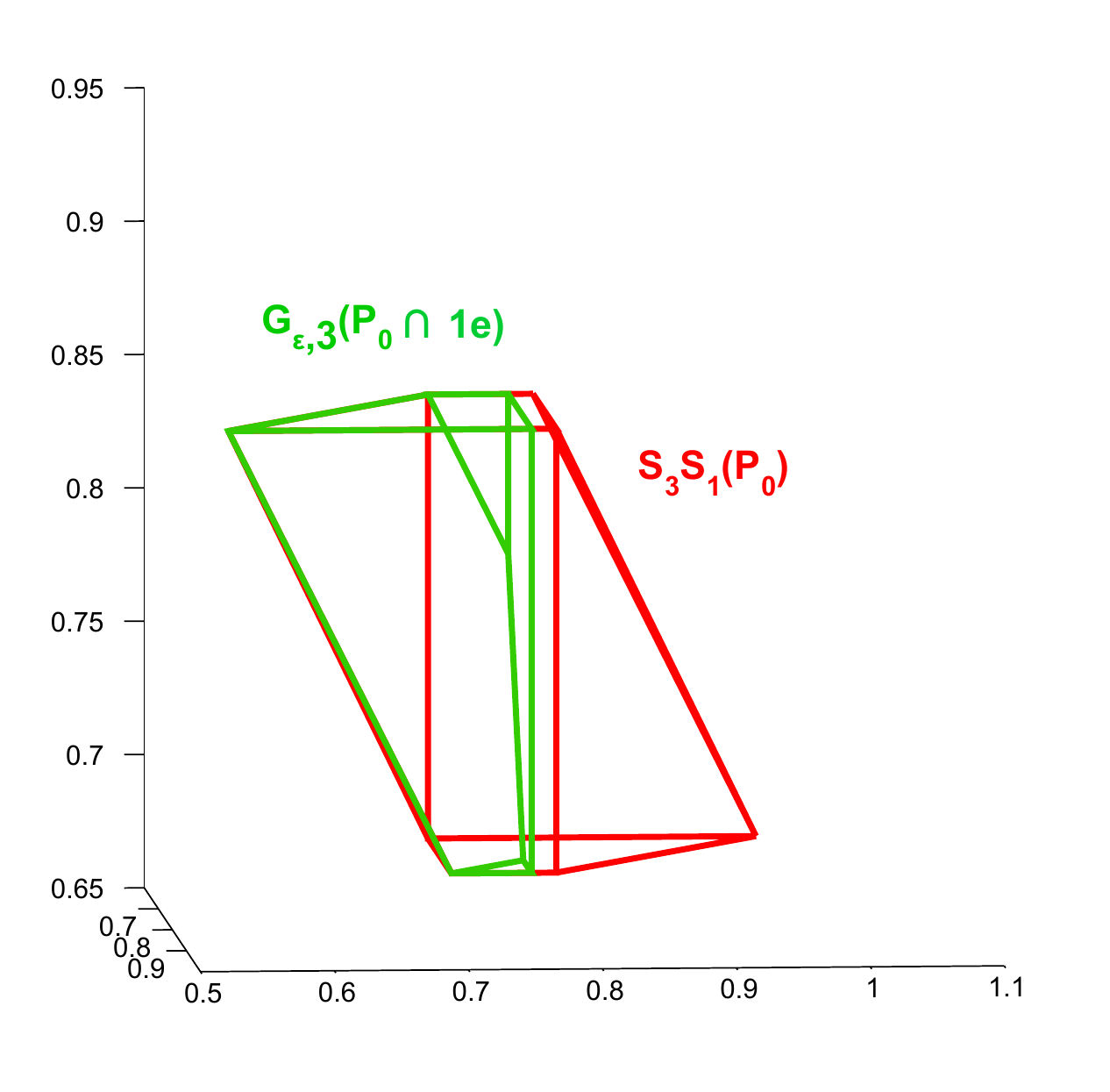}
   \caption{The image of $P_0 \cap 1e$. The case of the image of $P_0 \cap 8b$ is geometrically very similar.} \label{Fig1e8b}
   \end{subfigure}
  \caption{Images of $P_0 \cap 4b$ and $P_0 \cap 1e$ for $\varepsilon=0.32 > 1-\frac{\sqrt{2}}{2}$.} 
  \end{figure}

We also see that $G_{\varepsilon,3}(P_0 \cap 5b) \subseteq S_2(P_0)$ if and only if $L_{\varepsilon}^3(1-\varepsilon/2) \leq L_{\varepsilon}(1-\varepsilon/2)$, which gives the condition \eqref{sqrt2} on $\varepsilon$. For an illustration, see \Cref{Fig45b}.

Similarly, one can also see that $G_{\varepsilon,3}(P_0 \cap 1e) \subseteq S_3S_1(P_0)$ and $G_{\varepsilon,3}(P_0 \cap 8b) \subseteq S_4S_1(P_0)$ if and only if $L_{\varepsilon}^3(1-\varepsilon/2) \leq L_{\varepsilon}(1-\varepsilon/2)$. These give the condition \eqref{sqrt2} on $\varepsilon$ once again. For an illustration, see \Cref{Fig1e8b}.

%

\vspace{1cm}
\emph{E-mail address: } \tt selley@math.bme.hu

\end{document}